\documentclass[hidelinks,onefignum,onetabnum]{siamart220329}

\usepackage{lipsum}
\usepackage{amsfonts}
\usepackage{amsmath}
\usepackage{graphicx}
\usepackage{epstopdf}
\usepackage{todonotes}
\usepackage{algorithm}
\usepackage{algorithmic}
\usepackage{multirow}
\usepackage{booktabs}
\usepackage{subcaption}
\usepackage[utf8]{inputenc}
\usepackage{own_commands}

\ifpdf
\DeclareGraphicsExtensions{.eps,.pdf,.png,.jpg}
\else
\DeclareGraphicsExtensions{.eps}
\fi

\newsiamremark{remark}{Remark}
\newsiamremark{hypothesis}{Hypothesis}
\crefname{hypothesis}{Hypothesis}{Hypotheses}
\newsiamthm{claim}{Claim}

\headers{Low-rank computation for MOGPs}{S.~Esche and M.~Stoll}

\title{Low-rank computation of the posterior mean in Multi-Output Gaussian Processes\thanks{Submitted to the editors \today.
}}

\author{ Sebastian Esche\thanks{Chair of Scientific Computing, Department of Mathematics, University of Technology Chemnitz, 
		(\email{sebastian.esche@mathematik.tu-chemnitz.de}).}
	\and
	Martin Stoll\thanks{Chair of Scientific Computing, Department of Mathematics, University of Technology Chemnitz,  (\email{martin.stoll@mathematik.tu-chemnitz.de}).}}

\usepackage{amsopn}

\ifpdf
\hypersetup{
  pdftitle={Low-rank computation of the posterior mean in Multi-Output Gaussian Processes},
  pdfauthor={S.~Esche, M.~Stoll }
}
\fi

\begin{document}

\maketitle

\begin{abstract}
Gaussian processes (GP) are a versatile tool in machine learning and computational science. We here consider the case of multi-output Gaussian processes (MOGP) and  present low-rank approaches for efficiently computing the posterior mean of a MOGP. Starting from low-rank spatio-temporal data we consider a structured covariance function, assuming separability across space and time. This separability, in turn, gives a decomposition of the covariance matrix into a Kronecker product of individual covariance matrices. Incorporating the typical noise term to the model then requires the solution of a large-scale Stein equation for computing the posterior mean. For this, we propose efficient low-rank methods based on a combination of a \lrpcg\ method (\lrpcg) with the Sylvester equation solver \kpik\ adjusted for solving Stein equations. We test the developed method on real world street network graphs by using graph filters as covariance matrices. Moreover, we propose a degree-weighted average covariance matrix, which can be employed under specific assumptions to achieve more efficient convergence.
\end{abstract}

\begin{keywords}
Gaussian Process Regression, Multi-output Gaussian Process, Conjugate Gradient, Low-Rank, Stein equation, Sylvester equation
\end{keywords}

\begin{MSCcodes}
68Q25, 68R10, 68U05
\end{MSCcodes}

\section{Introduction}
Gaussian processes (GP) have received much attention in recent time due to their sophisticated prediction and uncertainty quantification capabilities (cf. \cite{duvenaudAutomaticModelConstruction2014, gramacySurrogatesGaussianProcess2020, rasmussenGaussianProcessesMachine2006, williamsGaussianProcessesRegression1995}). GPs are model driven due to its encoding mainly by the covariance function and mean function. As a result, GPs are therefore applied in a wide field of application, e.g. computational fluid dynamics or finite element analysis (see \cite{liuRemarksMultioutputGaussian2018} for an overview of applications).

Standard Gaussian processes are designed for scalar outputs but many application require vector-valued output. For applications with multiple outputs there might be complex relationships between the different output variables that can be addressed by using a Multi-Output Gaussian Process (MOGP), which transfers information across related outputs to improve prediction and uncertainty quality \cite{liuRemarksMultioutputGaussian2018, zhiGaussianProcessesGraphs2024}. More precisely, Gaussian processes (GPs) are commonly defined over index set $\mathbb{R}^n$, the index set can, in principle, be chosen arbitrarily. In particular, a GP defined over $\mathbb{R}^n \times\Wc$, where $\Wc$ is a set of output indices, naturally yields a MOGP. \cite{alvarezKernelsVectorValuedFunctions2012}.
Compared to a collection of single GPs a MOGP has the advantage of allowing cross covariance between output dimensions \cite{alvarezComputationallyEfficientConvolved2011}, which allows the incorporation of additional knowledge into the covariance model.

The outputs of a MOGP can also be interpreted as nodes of a graph \cite{dunsonGraphBasedGaussian2022} or are naturally posed on graphs \cite{rodriguez-denizUrbanNetworkTravel2017a}.  This structure offers a wide range of modeling capabilities, including the theory of predicting graph signals and spectrum based graph filtering \cite{liStochasticDeepGaussian2020, nguyenDetectingLowpassGraph2024, sandryhailaDiscreteSignalProcessing2013, smolaKernelsRegularizationGraphs2003, venkitaramanPredictingGraphSignals2019, venkitaramanGaussianProcessesGraphs2020, zhiGaussianProcessesGraphs2023}.

In this paper, we assume this MOGP setting with index set $\RE^n\times\Wc$, and graph interpretation of the outputs. We use therefore both a continuous dimension $\RE^n$ as well as a discrete dimension $\Wc$. For $n=1$ we have a spatio-temporal model, which e.g. can be used for describing quantities on a network evolving over time.
\begin{figure}
     \centering
     \caption{Perspectives in space and time for graph-based data.\label{fig:graphs}}
	\begin{subfigure}[b]{0.4\textwidth}
		\centering
		\includegraphics[height=5cm]{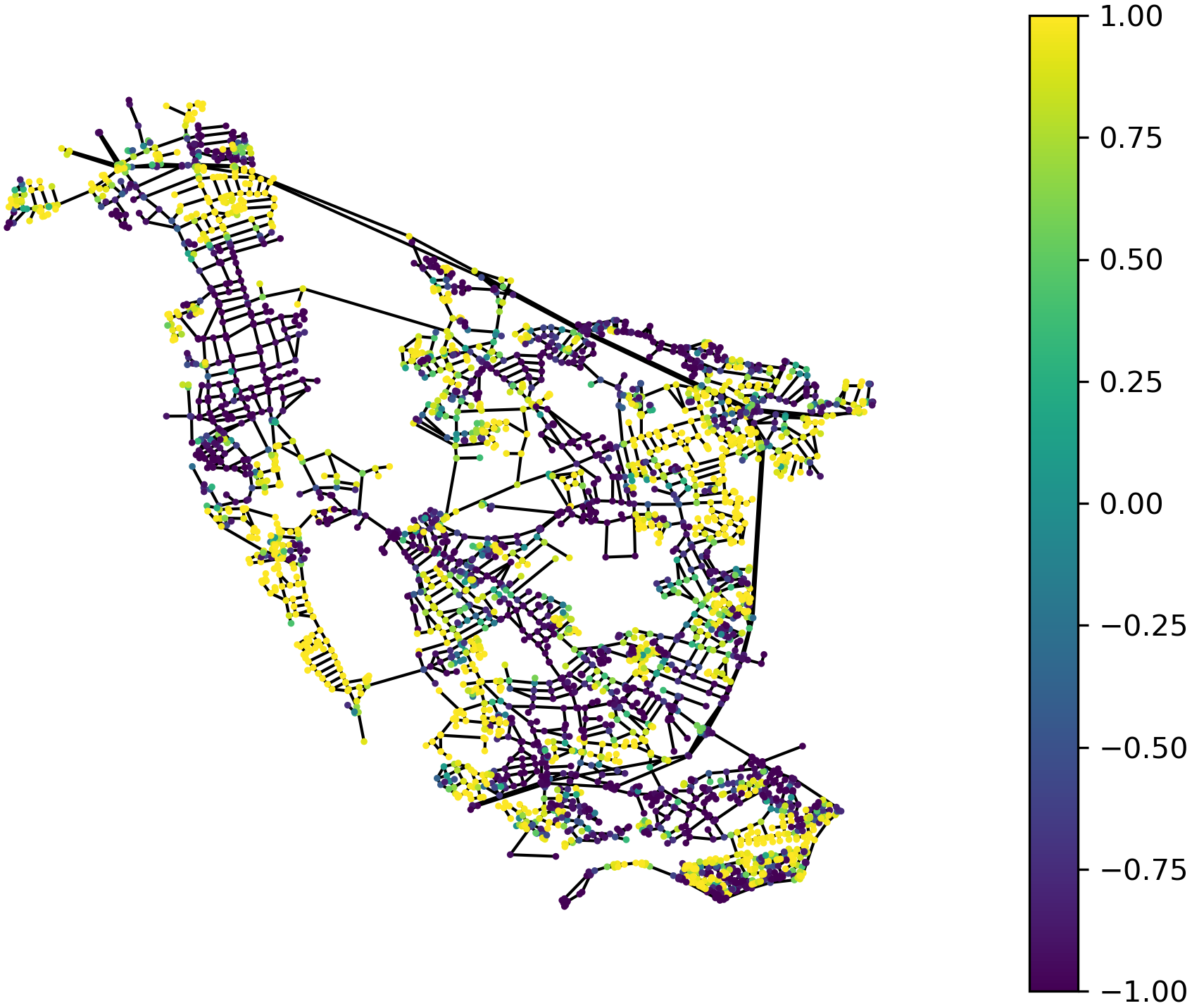}
		\caption{Signal on a street network graph}
		\label{fig:signaloxford}
	   \end{subfigure}
	\begin{subfigure}[b]{0.4\textwidth}
		\centering
		\includegraphics[height=5cm]{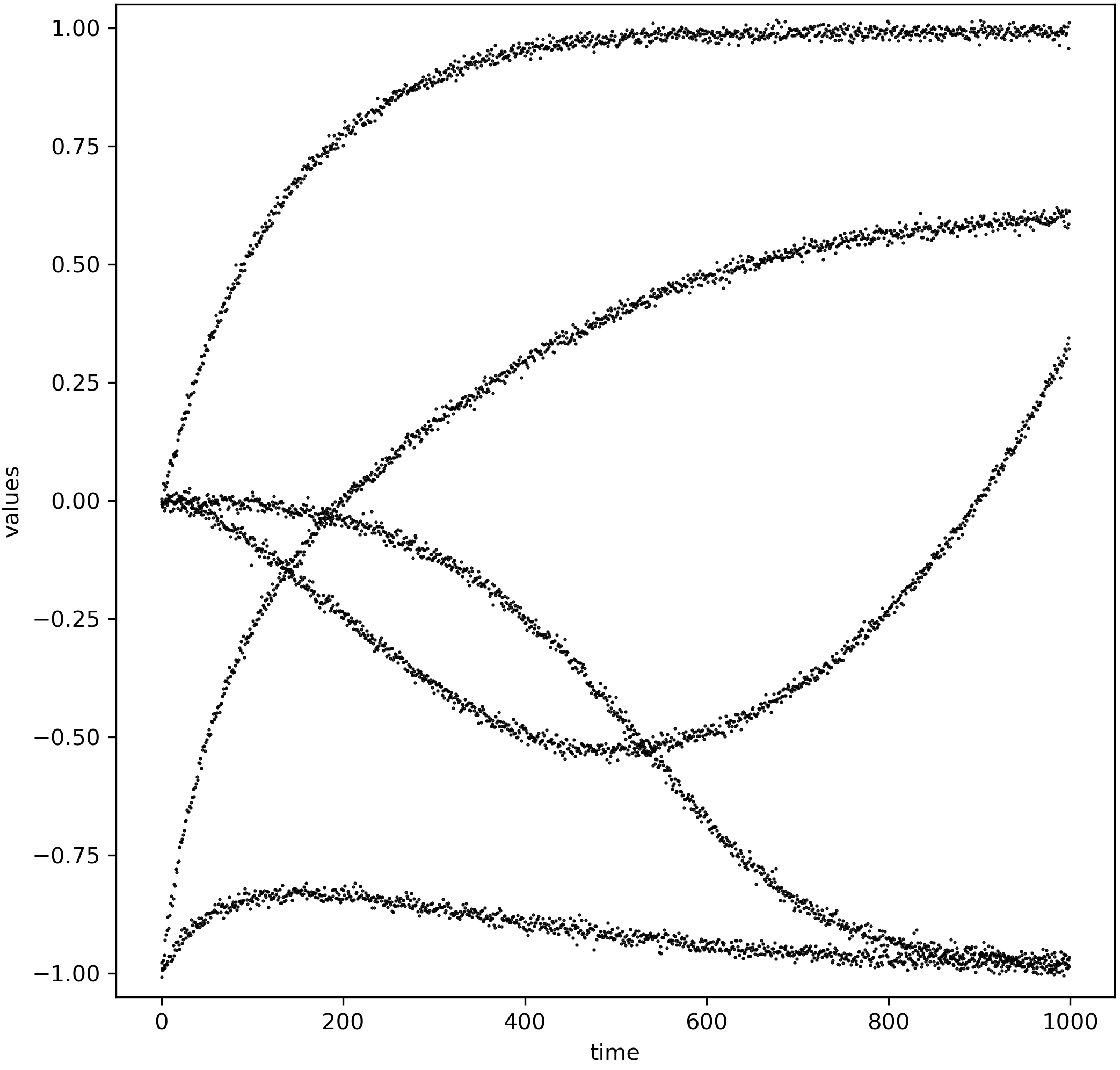}
		\caption{Multiple output values over time}
		\label{fig:timeseriesoxford}
   \end{subfigure}
\end{figure}
\Cref{fig:signaloxford,fig:timeseriesoxford} show simulated spatio-temporal data in different perspectives. \Cref{fig:signaloxford} shows the street network of Oxford with values at all nodes ranging from $-1$ to $1$ at a fixed time point. The collection of these values on nodes is also called a \textit{graph signal}. \Cref{fig:timeseriesoxford} shows evolution of the nodal values, i.e. a time series, for $5$ different randomly chosen nodes. Since we assume the index set to be $\RE^n\times\Wc$, for a finite subset of times it becomes a grid, for which the data can be stored in a matrix.

Large amounts of data occur in many real-world applications, e.g. due to a large amount of signals sampled or measured over time or the underlying graph contains millions of nodes, or both. Since the total size of the used data for the MOGP is the product of the number of signals and the number of nodes, the size can quickly become prohibitively large.

In general, an efficient approach to handle these sizes are low-rank methods \cite{bachmayrLowrankTensorMethods2023, breitenLowRankSolversFractional2016, grasedyckLiteratureSurveyLowrank2013, kressnerLowRankTensorKrylov2011, stollLowRankTimeApproach2015}, where the quantities of interest are approximated in a low-rank format. In these approaches the rank of the approximation typically depends on the properties of the solution vectors and the operators representing the covariance and correspondingly certain eigenvalue or singular value properties. We note that low-rank approximations of kernel matrices is not considered here as we explicitly focus on the utilization of the Kronecker structure. Approaches for low-rank approximations of the kernel matrix such as inducing point methods restrict the training data to a suitable subset-- the \textit{inducing points}. These methods are also known as sparse models; see, for example, \cite{burtConvergenceSparseVariational2020, nguyenCollaborativeMultioutputGaussian2014,  quinonero-candelaUnifyingViewSparse2005}. If we assume that the temporal covariance is independent of the covariance between outputs, the total covariance matrix can be expressed as the Kronecker product of the time covariance matrix $\KKI$ and the output covariance matrix $\KKO$ \cite{flaxmanFastKroneckerInference2015, nicksonBlitzkrigingKroneckerstructuredStochastic2015, stegleEfficientInferenceMatrixvariate2011}.
This Kronecker factorization reduces both storage requirements and computational costs.

If we incorporate the typical noise term $\sig^2$, together with the Kronecker structure it leads to the challenge of solving the system
\begin{align}
	\left(\KKI\otimes\KKO + \sig^2\II\right)\xx=\yy.
\end{align}
This is equivalent to solving a Stein equation. For modest sized matrices the Bartels-Stewart method \cite{bartelsAlgorithm432C21972} solves the Stein equation efficiently, but for large matrices this approach quickly becomes infeasible. In this case, we start with a low-rank right hand side and propose the use of an iteration method like the \kpik~method by Simoncini \cite{simonciniNewIterativeMethod2007} a state-of-the-art Krylov subspace solver. In this work, we study low-rank solvers, namely, an inexact \kpik\ solver and a \lrpcg method, which relies on an inexact \kpik~as preconditioner and controls the rank of the solution.

\subsection*{Outline}
The classical theory of GP regression is shown in \Cref{sec:gaussian_processes} and its generalization to MOGPs in \Cref{sec:mogp}. We draw a line between MOGPs and the related GMRFs and their numerical advantages in \Cref{sec:relation_to_gaussian_markov_random_fields}. In \Cref{sec:graph_filtering}, we show the connection of MOGPs to graph signal filtering and its use for covariance matrices. In addition, we propose a numerical beneficial covariance matrix, the Degree-Weighted Average, in \Cref{sec:degree_weighted_average}. For small sized problems the posterior mean of a MOGP can be computed by the eigendecomposition of the covariance matrices, which we illustrate in \Cref{sec:eigendecomposition}. One contribution of this paper is the derivation of the \lrpcg~method in \Cref{sec:iterative_solvers}. We demonstrate its capabilities and explain implementation details and the results in \Cref{sec:numerical_experiments}. See \Cref{sec:conclusion} for the conclusion.

\subsection*{Important concepts}
A crucial tool in the derivation of our method is the Kronecker product
\begin{align}
	\AA\otimes\BB = \begin{bmatrix}
		a_{11}\BB&\hdots&a_{1m}\BB\\
		\vdots&\ddots&\vdots\\
		a_{n1}\BB&\hdots&a_{nm}\BB
	\end{bmatrix}
\end{align}
which is very much intertwined with the $\mathrm{vec}$ operator defined via
\begin{align}
	\myvec{\AA} = \begin{bmatrix}
		\aa_{1}\\
		\vdots\\
		\aa_{m} \\
	\end{bmatrix}.
\end{align}
In order to derive certain equivalence relations and to motivate the low-rank approaches dicussed later we additionally state the following
\begin{align}\label{eqn:kronecker_vec}
	\left(\AA\otimes \BB\right)\myvec{\CC}=\myvec{\BB\CC\AA^\top}.
\end{align}
This simple relation becomes handy if the matrix $\CC$ is given in low-rank form or can be approximated well by a low-rank decomposition, i.e. $\CC\approx \UU_{\CC}\VV_{\CC}^\top$, as then we obtain
\begin{align}
	\left(\AA\otimes \BB\right)\myvec{\UU_{\CC}\VV_{\CC}^\top}=\myvec{\BB\UU_{\CC}\VV_{\CC}^\top \AA^\top}.
\end{align}
This requires the application of $\BB$ to $\UU_{\CC}$ and $\AA$ to $\VV_{\CC},$ so that $\BB\UU_{\CC}\VV_{\CC}^\top \AA^\top$ itself is the product of two low-rank matrices. This has been exploited heavily for designing low-rank solvers \cite{stollLowRankTimeApproach2015,kressnerLowRankTensorKrylov2011,grasedyckLiteratureSurveyLowrank2013,bachmayrLowrankTensorMethods2023,breitenLowRankSolversFractional2016}.  We will come back to this in the later derivation of our approach.

\subsection*{Notation}
Throughout this paper, we denote matrices by bold uppercase Roman letters, e.g., $\AA$, and vectors by bold lowercase Roman letters, e.g., $\aa$. In certain cases, we also use calligraphic letters, e.g., $\Bc$, to denote matrices when explicitly stated; otherwise, such letters refer to sets. An exception is made for input-output vectors and matrices, which we denote using unbold letters, such as $\Xu, \xu$ for inputs and $\Yu, \yu$ for outputs. This convention allows us to clearly distinguish them from unknown matrix variables like $\XX$ in matrix equations or from iterates like $\xx_k$ in iterative solvers. Furthermore, we distinguish between subscripts used for indexing multiple variables, such as $\xx_i$, and subscripts indicating components of a single vector, written as $(\xx)_i$. We use an asterisk to distinguish target data, e.g.  $\Xu^\ast, \Yu^\ast, \Vc^\ast$, from training data, e.g. $\Xu, \Yu, \Vc$. We denote the posterior by an asterisk $\MUp{}, \SIGMAp{}$ apart from the prior $\MU, \SIGMA$, which we will introduce in more detail later. The (column-wise) vectorization will be denoted by a $\mathrm{vec}$-subscript, e.g. $\yuv = \myvec{\Yu}$. The evaluation of a vector valued covariance function $k$ with matrix inputs $\AA=[\aa_1, \dots, \aa_n]\in\RE^{c\times n}, \BB=[\bb_1, \dots, \bb_n]\in\RE^{d\times m}$ is defined column wise by
\begin{align}\label{eqn:matrixcovariance}
	k(\AA,\BB) :=
	\begin{bmatrix}
		k(\aa_1,\bb_1) & \cdots & k(\aa_1,\bb_{m}) \\
		\vdots & & \vdots \\
		k(\aa_{n},\bb_1) & \cdots & k(\aa_{n},\bb_{m})
	\end{bmatrix},
\end{align}
or, in case of sets by the pairwise evaluation of their elements.

\section{Standard Gaussian Process Regression}\label{sec:gaussian_processes}
We here review some of the basic properties of Gaussian processes but refer to \cite{rasmussenGaussianProcessesMachine2006} for a detailed introduction. We start by assuming for $\Ic$ to be an arbitrary set. A family of random variables $(X_{t})_{t\in \Ic}$ with values in $\RE$ is called a \textit{Gaussian process (GP)} if for all $n\in \NN$ and for all $t_1, \dots , t_n\in \Ic$ it holds that $(X_{t_1}, \dots, X_{t_n})$ is $n$-dimensional normally distributed. Let $(X_t)_{t\in \Ic}$ be a GP with index set $\Ic$. The function $\Fc:\Ic\to \lbrace X_t : t\in \Ic \rbrace$ with $\Fc(t) := X_{t}$ is called function view of the GP. $\Fc$ can be seen as a random function, since each realization of $\Fc$ gives a (deterministic) function $f:\Ic\to\RE$. The function $m(t) := \EW[\Fc(t)]$ mapping $\Ic\to\RE$ is called \textit{mean function} of $\Fc$. The function $k(t, t') := \COV(\Fc(t), \Fc(t'))$ mapping $\Ic\times \Ic\to\RE$ is called \textit{covariance function} of $\Fc$ (cf. \cite{rasmussenGaussianProcessesMachine2006}). A GP is called \textit{centered} if it holds for its mean function $m(t) = 0$ for all $t\in \Ic$. The mean and covariance function uniquely determine the GP \cite{klenkeProbabilityTheoryComprehensive2014}.

The problem we are interested in is to efficiently compute the posterior distribution for unknown target values encoded in $\Fcp$ of a prior distribution $\Fc$ (both are GPs) in two steps. First, we determine $\Fc$ by choosing a mean- and a covariance function according to the model parameters. Here one can incorporate prior knowledge into the covariance function depending on the understanding of the generation of the data or possibly underlying physics. Second, we then infer $\Fcp$ from the values of $\Fc$ by a set of known training data. This is process is known as \textit{Gaussian process (GP) regression}, see \cite{williamsGaussianProcessesRegression1995}.

Let $k(\xu, \xu')$ be a chosen covariance function with $\xu, \xu'$ being elements of the index set $\Ic=\RE^c$. We assume $\Fc\sim\GP(0, k(\xu, \xu'))$ to be a centered GP prior and $\Fcp$ a GP posterior for which we want to infer the distribution. Thus, we assume for an arbitrary training input $\xu\in\RE^c$ that the corresponding output $\yu\in\RE$ is a realization of $\Fc$ and a noise term $\eps\sim\NV(0,\sig^2)$
\begin{align}\label{eqn:model}
	\yu = \Fc(\xu) + \eps.
\end{align}
Here, we assume $\eps$ to be i.i.d. with the variance $\sig^2>0$ to model uncertainty of the training data. Similarly, we assume for an arbitrary unknown target input $\xu^\ast\in\RE^c$ that the corresponding output is given as
\begin{align}
	\yu^{\ast} = \Fcp(\xu^\ast).
\end{align}

By definition of the prior distribution, we have $\COV(\Fc(\xu), \Fc(\xu'))=k(\xu,\xu')$. This means that for arbitrary inputs $\xu, \xu'$ and their corresponding outputs $\yu, \yu'$ using \Cref{eqn:model} and the independence of the noise gives
\begin{align}\label{eqn:covariances1D}
	\COV(\yu, \yu')=\COV(\Fc(\xu) + \eps, \Fc(\xu') + \eps')=k(\xu,\xu') + \delta_{\xu,\xu'} \sig^2.
\end{align}
Here, the error term appears only if $\xx = \xx'$, which is enforced by the Kronecker delta $\delta$, being one if $\xx = \xx'$ and zero otherwise. Consequently, this term vanishes entirely in the training-target and target-target covariance.

Let the \textit{training} data be given as $n$ input-output pairs $\lbrace\xu_i, \yu_i\rbrace_{i=1}^n$ and we denote the \textit{target} data as $n^\ast$ pairs $\lbrace\xu_j^\ast, \yu_j^\ast\rbrace_{j=1}^{n^\ast}$. We structure the data in the following way in matrix form as
\begin{equation}
	\begin{aligned}
		\Xu &= \begin{bmatrix}
			\vert & \vert & & \vert \\
			\xu_1 & \xu_2 & \dots & \xu_{n} \\
			\vert & \vert & & \vert \\
		\end{bmatrix} \in \RE^{c\times n}, \\
		\yu &= \begin{bmatrix}
			\yu_1 & \yu_2 & \dots & \yu_n
		\end{bmatrix}^{\top} \in \RE^n,
	\end{aligned}
	\quad
	\begin{aligned}
		\Xu^{\ast} &= \begin{bmatrix}
			\vert & \vert & & \vert \\
			\xu_1^\ast & \xu_2^\ast & \dots & \xu_{n^\ast}^\ast \\
			\vert & \vert & & \vert \\
		\end{bmatrix} \in \RE^{c\times n^\ast}, \\
		\yu^{\ast} &= \begin{bmatrix}
			\yu_1^\ast & \yu_2^\ast & \dots & \yu^\ast_{n^\ast}
		\end{bmatrix}^{\top}  \in \RE^{n^\ast}.
	\end{aligned}
\end{equation}
Here, we deviate from our matrix notation for input output matrices and using non bold characters, to distinguish it, from the later occurring matrix $\XX$ in the matrix equations. The training data $\Xu, \yu$ and the target inputs $\Xu^\ast$ are given where the target outputs $\yu^\ast$ are unknown. Namely, we want to infer their distribution. Usually, we interpret the $\xu$ input as a time, but for our experimental model in \Cref{sec:numerical_experiments} we also examine a case where it is reasonable to choose $\xu\in\mathbb{R}^c$ for $c>1$.

We define the following covariance matrices for the training and target inputs according to \Cref{eqn:matrixcovariance}
\begin{align}
	\KK:=k(\Xu,\Xu),\quad
	\KK^{\ast}:=k(\Xu,\Xu^\ast),\quad
	\KK^{\ast\ast}:=k(\Xu^\ast,\Xu^\ast).
\end{align}
The covariances of the components of $\yu$ and $\yu^\ast$ can compactly be denoted by	
\begin{align}
	\COV(\yu,\yu) = \KK + \sig^2\II,\quad
	\COV(\yu,\yu^\ast)=\KK^{\ast},\quad
	\COV(\yu^\ast,\yu^\ast)=\KK^{\ast\ast},
\end{align}
where the noise error $\sig^2\II$ emerges on the diagonal from~\Cref{eqn:covariances1D}. Using this we get the following prior distribution of the joint vector of $\yu$ and $\yu^\ast$:
\begin{align}\label{dist:noisy}
	\begin{bmatrix}
		\yu \\
		\yu^\ast
	\end{bmatrix}
	\sim \NV\left( 0,
	\begin{bmatrix}
		\KK + \sig^2\II & \KK^{\ast}\\
		{\KK^\ast}^\top & \KK^{\ast\ast}
	\end{bmatrix}
	\right).
\end{align}

For the prediction at the target inputs $\Xu^{\ast}$ we have to compute the posterior distribution, e.g. compute the mean vector and covariance matrix of $\yu^\ast$. To calculate the posterior distribution of $\yu^\ast$ with the knowledge of $\yu$ we receive by the conditional distribution of multivariate normal distribution
\begin{equation}\label{eqn:regression_noisy}
	\begin{aligned}
		\yu^\ast\vert\yu&\sim \NV(\MUp, \SIGMAp{}), \text{ with}
	\end{aligned}
	\begin{aligned}
		\MUp&={\KK^\ast}^\top[\KK+\sig^2\II]^{-1}\yu, \\
		\SIGMAp{}&=\KK^{\ast\ast}-{\KK^\ast}^\top[\KK+\sig^2\II]^{-1}\KK^\ast,
	\end{aligned}
\end{equation}
which is the mean vector and covariance matrix of the posterior, respectively.

\section{Multi-output Gaussian process regression}\label{sec:mogp}
A wide variety of problems contain relational information between objects which are generally modeled by graphs or do have a natural graph structure.

In the following, we will extend the previously defined Gaussian process model to the multi-output case, which are applicable to graphs, see \cite{liuRemarksMultioutputGaussian2018}. For this consider the following definition: Let $\Ic$ be an arbitrary set of input indices and let $\Wc=\lbrace 1,\dots, w\rbrace$ be a discrete set of output indices. A \textit{Multi-output Gaussian process} (MOGP) denoted by $\Fc(\xu, u)$ is a GP with index set $\tilde{\Ic}=\Ic\times \Wc$. So the output components are simply handled as an additional index of a GP. Hence, a MOGP is a GP that equips the random variables with additional output indices. We will see later that $\Wc$ could encode the nodes and $\kO$ the weights of the edges of a graph.

For $\Sc\subset\Wc$ we use subscript notation $\Fc_\Sc$ for expressing the vector with component $\Fc_s:=\Fc(\cdot,s)$ at position $s$ for all $s\in\Sc$. In particular $\Fc_\Wc:=\vek{\Fc_1 & \dots & \Fc_{w}}^\top$ is the vector of all output components. Furthermore, we distinguish between $m$ output components belonging to training data denoted by $\Vc\subset \Wc$  and $m^\ast$ output components belonging to target data denoted by $\Vc^\ast\subset \Wc$. It must hold $\Wc=\Vc\cup \Vc^\ast$ and because of possibly overlapping $\Vc$ and $\Vc^\ast$ in general $m+m^\ast\leq w$ holds true.

Let $\kI$ be a covariance function on the input space $\Ic=\RE^c$ and $\kO$ a covariance function on the discrete set of output dimensions $\Wc$. A covariance function $k$ is called \textit{separable} if it holds in general that
\begin{align}
	k((\xu, u),(\xu', u'))=\kI(\xu,\xu')\kO(u, u').
\end{align}
We call $\kI$ \textit{input covariance function} and $\kO$ \textit{output covariance function}.

If we want to do regression, similar to the previous section, we assume a centered MOGP prior on $\RE^{c}\times \Wc$
\begin{align}
	\Fc\sim\GP\left(0,k\left((\xu, u),(\xu', u')\right)\right)
\end{align}
with determined separable covariance function $k$. Therefore, it holds
\begin{equation}\label{eqn:mogp_covariance}
	\begin{aligned}
		\COV\left(\Fc(\xu, u), \Fc(\xu', u')\right) = k\left((\xu, u),(\xu', u')\right) = \kI(\xu, \xu')\kO(u,u').
	\end{aligned}
\end{equation}
We assume for $u\in \Vc$ that the $u$-th component of the training output $(\yu)_u$ is a realization of $\Fc_u$ and a noise error $\eps\sim\NV(0, \sig^2)$ (i.i.d.) with some $\sig^2>0$. We assume for $u^{\ast}\in \Vc^\ast$ that the $u^{\ast}$-th component of the target output $(\yu^\ast)_{u^\ast}$ is a realization of $\Fcp_{u^\ast}$. We denote this model in the following component-wise view as
\begin{align}\label{eqn:multi_output_model_1}
	(\yu)_u = \Fc_u(\xu) + \eps, \quad
	(\yu^\ast)_{u^\ast} = \Fcp_{u^\ast}(\xu^\ast).
\end{align}
This leads to a vector view given by
\begin{align}\label{eqn:multi_output_model_2}
	\yu = \Fc_\Vc(\xu) + \EPS,\quad
	\yu^\ast = \Fcp_{V^\ast}(\xu^\ast),
\end{align}
where $\EPS\in\RE^m$ and its components are i.i.d.

The covariance between training output components in \Cref{eqn:multi_output_model_1} is given by \Cref{eqn:mogp_covariance}, and it holds
\begin{align}\label{eqn:mogp_cov_output}
	\COV((\yu)_u,(\yu')_{u'}) &= \COV(\Fc_u(\xu) + \eps,  \Fc_u(\xu') + \eps') = \kI(\xu, \xu')\kO(u,u')+\delta_{\xu,\xu'}\delta_{u,u'}\sig^2.
\end{align}
Since the errors are assumed to be independent between output dimensions we need to multiply the error term by $\delta_{u,u'}$. We obtain the same result for covariance between training-target or target-target components by canceling the noise term in \Cref{eqn:mogp_cov_output}.

We assume $n$ given training data pairs $\lbrace\xu_{i},\yu_i\rbrace_{i=1}^n$ with input $\xu_i\in\RE^c$ and output $\yu_i\in\RE^m$ and we denote the target data by $\lbrace\xu_{j}^\ast,\yu_{j}^\ast\rbrace_{j=1}^{n^\ast}$, where $\xu_{j}^\ast\in\RE^c$ is a given input vector with corresponding output $\yu_{j}^\ast\in\RE^{m^\ast}$ which is unknown.
We structure the data in matrices
\begin{equation}\label{eqn:mogp_matrices}
	\begin{aligned}
		\Xu &= \begin{bmatrix}
			\vert & \vert & & \vert \\
			\xu_1 & \xu_2 & \dots & \xu_{n} \\
			\vert & \vert & & \vert \\
		\end{bmatrix} \in \RE^{c\times n}, \\
		\Yu &= \begin{bmatrix}
			\vert & \vert & & \vert \\
			\yu_1 & \yu_2 & \dots & \yu_n \\
			\vert & \vert & & \vert \\
		\end{bmatrix} \in \RE^{m\times n},
	\end{aligned}
	\quad
	\begin{aligned}
		\Xu^{\ast} &= \begin{bmatrix}
			\vert & \vert & & \vert \\
			\xu_1^\ast & \xu_2^\ast & \dots & \xu_{n^\ast}^\ast \\
			\vert & \vert & & \vert \\
		\end{bmatrix} \in \RE^{c\times n^\ast}, \\
		\Yu^{\ast} &= \begin{bmatrix}
			\vert & \vert & & \vert \\
			\yu_1^\ast & \yu_2^\ast & \dots & \yu^\ast_{n^\ast}\\
			\vert & \vert & & \vert \\
		\end{bmatrix} \in \RE^{m^\ast\times n^\ast}.
	\end{aligned}
\end{equation}

Therefore, we distinguish between inputs $\Xu$ and outputs $\Yu$, where the training and target outputs  $\Yu$ and $\Yu^\ast$ can be on different sets of nodes $\Vc$ and $\Vc^\ast$.
We use the following abbreviations
\begin{equation}\label{eqn:mogp_KKI}
	\begin{aligned}
		\KKI = \kI(\Xu, \Xu),\quad \KKI^\ast = \kI(\Xu, \Xu^\ast),\quad \KKI^{\ast\ast} = \kI(\Xu^\ast, \Xu^\ast) \\
		\KKO = \kO(\Vc, \Vc),\quad \KKO^\ast = \kO(\Vc, \Vc^\ast),\quad \KKO^{\ast\ast} = \kO(\Vc^\ast, \Vc^\ast),
	\end{aligned}
\end{equation}
with $\Ic$ referring to input and $\Oc$ to output. We again store all training and target output components from the model in a single column vector by using the column-wise vectorization $\yuv = \myvec{\Yu}$ and $\yuva = \myvec{\Yu^{\ast}}$. This constructed vector is distributed according to
\begin{align}
	\vek{\yuv \\ \yuva} \sim \NV(0, \SIGMA),
\end{align}
where the covariance matrix $\SIGMA$ is given by
\begin{align}\label{eqn:cov_mogp}
	\SIGMA = \begin{bmatrix}
		\KKI\otimes\KKO & \KKI^\ast\otimes\KKO^\ast \\
		{\KKI^\ast}^\top\otimes{\KKO^\ast}^\top & \KKI^{\ast\ast}\otimes\KKO^{\ast\ast}
	\end{bmatrix}
	+ \begin{bmatrix}
		\sig^2\III\otimes\IIO & 0 \\
		0 & 0
	\end{bmatrix}.
\end{align}

To see that this equation is correct we break it down step by step and start with the definition of $\SIGMA$
\begin{align}
	\SIGMA = \begin{bmatrix}
		\COV(\yuv,\yuv) & \COV(\yuv,\yuva) \\
		\COV(\yuv,\yuva)^\top & \COV(\yuva,\yuva)
	\end{bmatrix}.
\end{align}
The four entries are matrix blocks given by the pairings of training and target variables. These blocks itself contain covariance matrices of the following form
\begin{align}
	\COV(\yuv,\yuv) = \begin{bmatrix}
		\COV(\yu_1,\yu_1) & \dots & \COV(\yu_1,\yu_n) \\
		\vdots & & \vdots \\
		\COV(\yu_n,\yu_1) & \dots & \COV(\yu_n,\yu_n)
	\end{bmatrix}.
\end{align}
Each of these pairs is again a covariance matrix between the components of the first vector against the components of the second vector. By our assumption each entry can be written as a product of two covariance functions, e.g.
\begin{equation}
	\begin{aligned}\label{eqn:output_vector_cov}
		\COV(\yu,\yu') &=
		\begin{bmatrix}
			\kI(\xu,\xu')\kO(1,1) & \dots & \kI(\xu,\xu')\kO(1,m) \\
			\vdots & & \vdots \\
			\kI(\xu,\xu')\kO(m,1) & \dots & \kI(\xu,\xu')\kO(m,m)
		\end{bmatrix} +
		\begin{bmatrix}
			\sig^2 &  & 0 \\
			 & \ddots &  \\
			0 &  & \sig^2
		\end{bmatrix},
	\end{aligned}
\end{equation}
which can be simplified to
\begin{align}\label{eqn:output_vector_cov_2}
	\COV(\yu,\yu')= \kI(\xu,\xu')\KKO + \sig^2\IIO.
\end{align}
From this we can conclude for the first block in $\SIGMA$
\begin{align}\label{eqn:output_vector_cov_3}
	\COV(\yuv,\yuv) = \KKI\otimes\KKO + \sig^2\III\otimes\IIO,
\end{align}
and we can proceed in the same way for the other blocks. Because we considered training-training covariance here, the noise term emerged, which is not the case for training-target or target-target covariance. Note, that it holds $\III\otimes\IIO=\II_{mn}$ but we will later use the Kronecker property if convenient.

For the posterior distribution we infer the training data by using the conditional distribution of the multivariate normal distribution
\begin{align}
	\yuva\vert\yuv \sim \NV(\MUp, \SIGMAp{})
\end{align}
for which we require the computation of
\begin{align}
	\MUp &= \left( {\KKI^\ast}^\top\otimes{\KKO^\ast}^\top\right) \lbrack \KKI\otimes\KKO + \sig^2\III\otimes\IIO \rbrack^{-1}\yuv, \\
	\SIGMAp{} &= \KKI^{\ast\ast}\otimes\KKO^{\ast\ast} - \left({\KKI^\ast}^\top\otimes{\KKO^\ast}^\top\right) \lbrack \KKI\otimes\KKO + \sig^2\III\otimes\IIO \rbrack^{-1} (\KKI^\ast\otimes\KKO^\ast).
\end{align}

One of the key challenges in prediction of the posterior data points comes from solving systems with the matrix
\begin{align}\label{eqn:K}
	\Kc := \KKI\otimes\KKO + \sig^2\III\otimes\IIO,
\end{align}
which we will focus on in \Cref{sec:eigendecomposition,sec:iterative_solvers} of this paper. Note that the main computational cost in computing the posterior mean comes from solving the system
\begin{equation}
	\label{eqn:stein1}
	\Kc\myvec{\Xb}=\myvec{\CC} \Leftrightarrow  \KKO\Xb\KKI + \sig^2\XX=\CC
\end{equation}
using the symmetry of the kernel matrices. Here, $\CC$ encodes the training data of the MOGP. Solving these systems does not only occur when computing the mean but also in the optimization of the hyperparameters of the Gaussian processes \cite{williamsGaussianProcessesRegression1995}.

For the kernel on the input space $\KKI$ we use the well known squared exponential kernel
\begin{align}\label{eqn:se_kernel}
	\kO(\xu,\xu^\ast)=\sig_{w}^2\exp\left(-\frac{1}{2\ell^2}\norm{\xu-\xu^\ast}_2^2\right),
\end{align}
depending on the length hyperparameter $\ell$ and the variance hyperparameter $\sig_{w}$. While $\ell$ determines how quickly the covariance between inputs decays with distance, $\sig_{w}$ scales the overall amplitude of the function.

For the kernel on the outputs we want to incorporate graph information. One approach involves Gaussian Markov Random Fields (GMRFs), which rely on sparse precision matrices derived from graph structure. These models typically require direct access to the inverse of the covariance matrix. Although this is not available in our setting, GMRFs still offer valuable insights. To clarify their connection to graph signal filtering, which we use for the output kernel, we briefly introduce GMRFs next.

\section{Relation to Gaussian Markov Random Fields}\label{sec:relation_to_gaussian_markov_random_fields}
The theory of Gaussian Markov Random Fields (GMRF) \cite{rueGaussianMarkovRandom2005} is strongly connected to Gaussian Processes. Consider the so-called \textit{precision matrix} $\QQ$, which is defined as the inverse of the covariance matrix, i.e., $\QQ = \SIGMA^{-1}$. The precision matrix has the property of encoding the pairwise conditional dependencies: the corresponding entry of two variables in $\QQ$ is zero iff they are independent given all other variables. In fact, Gaussian Markov Random Fields (GMRFs) assume that the multivariate normal distribution is specified by the precision matrix. For the graph $\Gc=(\Wc,\Ec)$ of a GMRF the Markov property holds, which states that for entries of $\QQ$ holds
\begin{align}\label{eqn:markovproperty}
	\qq_{ij}\neq 0 \quad \Longleftrightarrow \quad \lbrace i,j\rbrace\in\Ec.
\end{align}
This means that nonzero entries in the precision matrix correspond precisely to the edges of the graph.

For the following Gaussian prior 
\begin{align}
	\vek{\yuv \\ \yuv^\ast}\sim\NV\left(0, \vek{\QQ_{11} & \QQ_{12} \\ \QQ_{12}^\top & \QQ_{22}}^{-1}\right),
\end{align}
the computation of the posterior mean becomes
\begin{align}
	\MUp = -\QQ_{22}^{-1}\QQ_{12}^\top\yuv.
\end{align}
By comparing this with \Cref{eqn:regression_noisy} where the covariance matrix of training points $\SIGMA_{11}$ has to be inverted, here the precision matrix of the targets $\QQ_{22}$ has to be inverted.

If the set of edges $\Ec$ of a GMRF is of modest size, this results in a sparse precision matrix, where the corresponding covariance matrix is dense in general. Therefore the ``big $n$ problem'', i.e. computational cost of $O(n^3)$ for solving with the covariance matrix, is becoming manageable using the precision matrix instead. If there are no adjacent target nodes at all, $\QQ_{22}$ becomes the diagonal matrix of node degrees and is therefore easy to invert.

If we model by using the covariance matrix $\SIGMA$ in \Cref{eqn:cov_mogp} we have in general no direct access to $\QQ$, i.e. we cannot take computational advantage from it. Constructing $\QQ$ from sketch is restricted to being symmetric and positive semi-definite, since it is the inverse of a symmetric and positive semi-definite covariance matrix. I.e. we cannot use an arbitrary matrix which fulfills \Cref{eqn:markovproperty}. This makes it hard, to construct precision matrices in general. An approach given in \cite{lindgrenExplicitLinkGaussian2011} shows, that there is an explicit link between GPs and GMRFs using a weak solution of a stochastic partial differential equation on an arbitrary triangulation. There, preprocessing is required, such as setting up the model for a triangulation and the GMRF representation. This approach is limited to a subset of the Matérn class covariance function, where the graph must be a triangulation. For more details see also \cite{borovitskiyMaternGaussianProcesses2021}.

Even though graph signal filtering in the following section does not make use of the Markov property like GMRFs do, both approaches rely on a sparse inverse covariance matrix, derived from a graph.

\section{Graph Signal Filtering and MOGPs}\label{sec:graph_filtering}
In order to motivate the choice of $\KKO$ we now introduce a graph signal approach given in \cite{zhiGaussianProcessesGraphs2023}, for other models using MOGPs see also \cite{liuRemarksMultioutputGaussian2018, stankovicIntroductionGraphSignal2019, venkitaramanGaussianProcessesGraphs2020}.
The graph filtering model is an important special case of the previously presented general MOGP model as we will see.

We use the setting given in \cite{zhiGaussianProcessesGraphs2023}. Let $\Gc=(\Wc,\Ec)$ be an undirected graph with $m$ vertices $\Wc$ and edges $\Ec$. We assume for the graph filtering that training and target output space are equal, it then holds that $\Vc=\Vc^\ast=\Wc$ and $m=m^\ast=w$. For a given adjacency matrix (or weight matrix) $\WW$ and diagonal degree matrix $\DD$ (or matrix of row sums of $\WW$) the graph Laplacian of $\Gc$ is defined by $\LL=\DD-\WW$. Let $\LL=\UU\LAMBDA \UU^\top$ be the eigendecomposition of the Laplacian where $\UU$ contains the eigenvectors and $\LAMBDA$ is the diagonal matrix of eigenvalues. We now consider data on $\Gc$, i.e. for each input $\xu\in\RE^c$ we get an unfiltered output $\ff\in\RE^m$ where we identify each output component of $\ff$ with a node of $\Gc$. In this context $\ff$ is also called a \textit{signal} (cf. \cite{shumanEmergingFieldSignal2013}). The graph Fourier transform of the signal is defined by $\UU^\top\ff$ and computes the spectrum of $\ff$. Let $g(\LAMBDA)$ be a non-negative function. The filtered signal is then given by
\begin{align}
	\yu=g(\LL)\ff=\UU g(\LAMBDA)\UU^\top\ff
\end{align}
and we now discuss how to select a suitable filter function.

The filtering model presented in \cite{zhiGaussianProcessesGraphs2023} uses a different approach of introducing covariance in multi-output Gaussian processes but leads to same results. We show, that both approaches are indeed equivalent. According to the filtering model, we call a MOGP $\Fc$ with output indices $\Wc$ a \textit{simple} MOGP if all of its output processes $\Fc_1, \dots, \Fc_m$ are themselves independent GPs with identical covariance function. Let $\kI$ be a covariance function of inputs from $\RE^c$ and let $\Fc\sim\GP(0, \kI(\xu,\xu')\delta_{u,u'})$ be a \textit{simple} MOGP, where $u, u'$ are nodes in $\Wc$. The covariance function $\kI(\xu,\xu')\delta_{u,u'}$ is separable and $\delta_{u,u'}$ ensures the independence of the output components. This leads to the following covariance matrix of the vector view
\begin{align}
	\COV(\Fc_\Wc(\xu), \Fc_\Wc(\xu')) = \kI(\xu,\xu')\II_m.
\end{align}

We now consider the matrix $\Bc\in\RE^{m\times m}$ called \textit{filtering matrix}. The covariance between outputs is then introduced by considering the filtered signal
\begin{align}\label{eqn:filtering_model}
	\yu = \Bc\Fc_\Wc(\xu)+\EPS.
\end{align}
The multiplication of $\Fc_\Wc$ by $\Bc$ leads to covariances within the components of $\yu$. Hence, the filtering by $g$ is obtained by multiplication by $\Bc$ where $\Fc_\Wc$ takes the role of $\ff$. It holds
\begin{align}
	\COV(\yu, \yu')=\COV(\Bc\Fc_\Wc(\xu)+\EPS, \Bc\Fc_\Wc(\xu')+\EPS') =\kI(\xu,\xu')\BcT +\delta_{\xu,\xu'}\sig^2\II_m.
\end{align}	
The structure of this system is the same as the one given in \Cref{eqn:output_vector_cov_2} and as a result this approach is equivalent to the model given in the previous section. If we let $\Bc=\vek{\bb_1 & \cdots & \bb_m}\in\RE^{m\times m}$ and $u,u'\in\Wc$ we can even give the explicit link between both models by
\begin{align}
	(\KKO)_{u, u'} =\kO(u, u')=\bb_u^\top \bb_{u'}.
\end{align}
In a nutshell, the following approaches are equivalent:
\begin{equation}
	\begin{aligned}
		\yu&=\Bc\Fc_\Wc(\xu)+\EPS \text{ with } \Fc\sim\GP\left(0, \kI(\xu, \xu')\delta_{u,u'}\right), \\
		\yu&=\Fc_\Wc(\xu)+\EPS \text{ with } \Fc\sim\GP\left(0, \kI(\xu, \xu')\bb_u \bb_{u'}^\top\right).
	\end{aligned}
\end{equation}

Both models lead to the overall system
\begin{align}
	\COV(\yuv, \yuv)=\KKI\otimes \BcT + \sig_{\eps}^2\III\otimes\II_m,
\end{align}
which is the same as \Cref{eqn:output_vector_cov_3} with $\BcT=\KKO$. As we assumed that training and target outputs are equal and therefore $\KKO=\KKO^\ast=\KKO^{\ast\ast}$, the posterior covariance matrix is
\begin{align}
	\SIGMA = \begin{bmatrix}
		\KKI\otimes \BcT & \KKI^\ast\otimes \BcT \\
		{\KKI^\ast}^\top\otimes \BcT & \KKI^{\ast\ast}\otimes \BcT
	\end{bmatrix}
	+ \begin{bmatrix}
		\sig_{\eps}^2\III\otimes\II_m & 0 \\
		0 & 0
	\end{bmatrix}.
\end{align}

As explained in \cite{venkitaramanGaussianProcessesGraphs2020} we can motivate the choice of the kernel on the output space $\KKO=\BcT$ by a smoothness consideration. The smoothness of a signal $\ell(\yu)=(\yu^\top\LL\yu)/(\yu^\top\yu)$ quantifies the similarity between the values on the nodes. If we want to filter high frequencies, which corresponds to smoothening the signal, we have to solve the problem
\begin{align}
	\yu &= \argmin_{\zz}\norm{\ff-\zz}^2+\alpha\zz^\top\LL\zz = (\II+\alpha\LL)^{-1}\ff.
\end{align}
The regularization parameter $\alpha\geq0$ weights the penalty of the smoothness. This motivates the choice $\Bc=(\II+\alpha\LL)^{-1}$ which is also called \textit{global filter}.

Other choices for graph filters are given by \cite{venkitaramanGaussianProcessesGraphs2020}, e.g.
\begin{itemize}
	\item Standard GP: $\Bc=\II$,
	\item Global filtering: $\Bc=(\II+\alpha \LL)^{-1}$,
	\item Local averaging: $\Bc=(\II+\alpha \DD)^{-1}(\II+\alpha \WW)$,
	\item Graph Laplacian regularization (pseudo-inverse): $\BcT=\LL^{+}$,
	\item Regularized Laplacian: $\BcT=(\II+\alpha \LLt)^{-1}$.
\end{itemize}
The normalized graph Laplacian is given by $\LLt=\DD^{-\frac{1}{2}}\LL\DD^{-\frac{1}{2}}$, which has eigenvalues only in the interval $[0,2]$.  Even without knowledge about eigenvalues of the choice of $\Bc$, the eigenvalues of $\BcT$ have to be greater or equal to zero. Therefore, the requirement of being positive semi-definite to serve as a covariance matrix is fulfilled. See also \cite{smolaKernelsRegularizationGraphs2003} for details.

\section{Degree-Weighted Average}\label{sec:degree_weighted_average}
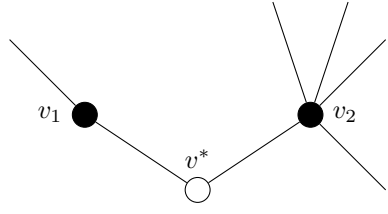
\begin{figure}
	\centering
	\begin{tikzpicture}
		\node[draw, circle, label=above:$v^\ast$] (C) at (0,0) {};
		\node[draw, circle, fill=black, label=left:$v_1$] (N1) at (-1.5,1) {};
		\node[draw, circle, fill=black, label=right:$v_2$] (N2) at (1.5,1) {};
		\draw (C) -- (N1);
		\draw (C) -- (N2);
		\draw (N1) -- ++(-1,1);
		\draw (N2) -- ++(1,1);
		\draw (N2) -- ++(1,-1);
		\draw (N2) -- ++(0.5,1.5);
		\draw (N2) -- ++(-0.5,1.5);
	\end{tikzpicture}
	\caption{Degree-Weighted Average: $v_1$ distributes $\frac{1}{2}$ and $v_2$ distributes $\frac{1}{5}$ of an abstract quantity to the (unknown) $v^\ast$. The proportions are assumed to be equal to all neighbors and therefore $1/\text{degree}$.}
	\label{fig:degreeweightedaverage}
\end{figure}

We want to propose another covariance matrix which can be applied for a special structure of the data but leads to an efficient way for solving with the matrix $\KKO$ and therefore enhances the computation of the posterior mean. In this case, we assume that the given data $\yu$ at training nodes $\Vc$ and the posterior mean $\MUp$ at target nodes $\Vc^\ast$ satisfy the stationary equation. This means that a quantity $y_j$ sitting in node $v_j$ distributes equally to all of its neighbors. From the perspective of such a neighbor $v_i$, the received quantity is $y_j/\deg(v_j)$. If we sum this up for all neighbors of $v_i$, we obtain the following relation
\begin{align}\label{eqn:weightedaverage}
	\matrix{\mu}^\ast_i = \sum_{j=1}^b \frac{\WW_{i j}}{\deg(v_j)}y_j,
\end{align}
where $\WW_{ij}$ encodes the adjacency between nodes $i$ and $j$, see \Cref{fig:degreeweightedaverage}. Note that this requires knowledge of all neighbors of the target nodes. If we cannot meet this requirement because some neighbors are target nodes itself, we get the more general system
\begin{align}\label{eqn:weightedaverage2}
	\matrix{\mu}^\ast_i = \sum_{v_j\in V} \frac{(\WW)_{i j}}{\deg(v_j)}y_j + \sum_{v_j\in V^\ast}\frac{(\WW)_{i j}}{\deg(v_j)}\mu_j^\ast,
\end{align}
which we have to solve for $\MUp$.

This linear equation system is closely connected to the stationary distribution
\begin{align}\label{eqn:stationary}
	\pi\PP=\pi\quad\Leftrightarrow\quad\pi(\PP-\II)=0,
\end{align}
where $\PP=\DD^{-1}\WW$ is the so called transition matrix and $\pi$ is a row vector corresponding to the nodes. For our purposes, we consider $\pi$ not as a distribution, i.e. we don't restrict $\pi\boldone=1$, but we define it as $\pi^\top:=\vek{\yu \\ \MUp}$. We will see that \Cref{eqn:weightedaverage2} is a consequence of the stationary distribution \Cref{eqn:stationary}. If we transpose \Cref{eqn:stationary} and insert the definitions, we obtain
\begin{align}
	\left(\WW\DD^{-1}-\II\right)\vek{\yu \\ \MUp}=0.
\end{align}
If we split the matrices according to the sizes of $\yu$ and $\MUp$, i.e.
\begin{align}
	\WW=\begin{bmatrix}
		\WW_{11} & \WW_{12} \\
		\WW_{12}^\top & \WW_{22}
	\end{bmatrix}, \quad
	\DD = \begin{bmatrix}
		\DD_1 & 0 \\
		0 & \DD_2
	\end{bmatrix},
\end{align}
we get the two equation systems
\begin{align}
	\left(\WW_{11}\DD_1^{-1}-\II\right)\yu + \WW_{12}\DD_2^{-1}\MUp &= 0 \\
	\WW_{21}\DD_1^{-1}\yu + \left(\WW_{22}\DD_2^{-1}-\II\right)\MUp &= 0,
\end{align}
where the identity matrices are of fitting sizes. As $\yu$ contains known variables the system is overdetermined. Hence, we drop the first and just use the second equation system. We solve this for $\MUp$
\begin{align}\label{eqn:mu_formula}
	\MUp &= \left(\II-\WW_{22}\DD_2^{-1}\right)^{-1}\WW_{21}\DD_1^{-1}\yu.
\end{align}
This is the matrix version of \Cref{eqn:weightedaverage2}. Pay attention that $\DD_2$ is using the sum of rows of $\vek{\WW_{21} & \WW_{22}}$ and not just of $\WW_{22}$. As we assume to have connections between $\Vc$ and $\Vc^\ast$ it holds $\WW_{21}\neq 0$. Therefore $\II-\WW_{22}\DD_2^{-1} = (\DD_2-\WW_{22})\DD_2^{-1}$ is not a product of a non singular Laplace factor.

We want to bridge to Gaussian process regression by considering a centered Gaussian Process. In this context the posterior mean is computed by $\MUp=\SIGMA_{12}^\top\SIGMA_{11}^{-1}\yu$. If we choose $\SIGMA_{11}=\DD_1$ and $\SIGMA_{12}=\WW_{12}\left(\II-\DD_2^{-1}\WW_{22}\right)^{-1}$ the posterior mean coincides with \Cref{eqn:mu_formula}. Therefore we achieve in the regression context the degree-weighted average from \Cref{eqn:weightedaverage2} by defining the covariance matrix as
\begin{align}\label{eqn:stat_cov_matrix}
	\SIGMA=\begin{bmatrix}
		\SIGMA_{11} & \SIGMA_{12} \\
		\SIGMA_{12}^\top & \SIGMA_{22}
	\end{bmatrix} &= \begin{bmatrix}
		\DD_1 & \MM \\
		\MM^\top & \SIGMAp{22} + \MM^\top\DD_1^{-1}\MM
	\end{bmatrix},
\end{align}
with $\MM=\WW_{12}\left(\II-\DD_2^{-1}\WW_{22}\right)^{-1}$. We treat the posterior covariance $\SIGMAp{22}$ as a free parameter which has to be symmetric and positive semi-definite. If it is positive definite it has a Cholesky decomposition $\SIGMAp{22}=\tilde{\LL}\tilde{\LL}^\top$. This gives the Cholesky decomposition of $\SIGMA=\LL\LL^\top$, which is therefore also positive definite
\begin{align}
	L=\begin{bmatrix}
		\DD_1^{\frac{1}{2}} & 0 \\
		\MM^\top\DD_1^{-\frac{1}{2}} & \tilde{\LL}
	\end{bmatrix}.
\end{align}
A similar reasoning can be applied when $\SIGMAp{22}$ is positive semi-definite. Since $\SIGMA$ is symmetric and positive semi-definite it is a covariance matrix.

For the case that there are no connected couples in $\Vc^\ast$ it holds $\WW_{22}=0$. If we additionally assume $\SIGMAp{22}=\DD_2-\WW_{21}\DD_1^{-1}\WW_{12}$ \Cref{eqn:stat_cov_matrix} reduces to
\begin{align}\label{eqn:simple_case}
	\SIGMA=\begin{bmatrix}
		\SIGMA_{11} & \SIGMA_{12} \\
		\SIGMA_{12}^\top & \SIGMA_{22}
	\end{bmatrix} &= \begin{bmatrix}
		\DD_1 & \WW_{12} \\
		\WW_{21} & \DD_2
	\end{bmatrix}.
\end{align}
and the posterior mean is computed by 
\begin{align}
	\MUp &= \WW_{21}\DD_1^{-1}\yu.
\end{align}
For this matrix we show positive semi-definiteness: The Gershgorin circle $i$ of $\SIGMA$ is given by its center at $\deg(v_i)$ and its radius
\begin{align}
	R_i=
	\left\{
	\begin{array}{lr}
		\sum_{j\in V^\ast}(\WW)_{i,j}, & \text{if } i \in V \\
		\sum_{j\in V}(\WW)_{i,j}, & \text{if } i \in V^\ast
	\end{array}
	\right\} \leq \sum_{j\in V^\ast}(\WW)_{i,j} + \sum_{j\in V}(\WW)_{i,j} = \deg(v_i). 
\end{align}
The Gershgorin circle theorem states that each eigenvalue $\lambda_i$ of $\SIGMA$ lies within its Gershgorin circle. Hence it holds $\lambda_i\in[0, 2\deg(v_i)],$ which implies $\lambda_i\geq0$ and so $\SIGMA$ is positive semi-definite.

The big advantage of this approach is that $\SIGMA_{11}=\DD_1$ is diagonal. This enables very fast multiplication and solving by $\KKO$ in the Stein equation, which leads to a fast iteration step for the solvers described in \Cref{sec:iterative_solvers}. Altough $\left(\II-\WW_{22}\DD_2^{-1}\right)^{-1}$ involves solving, it has to be done only once and just for a vector after the Stein equation is solved.

\section{Using the eigendecomposition}\label{sec:eigendecomposition}
When solving with the matrix 
\begin{align}\label{eqn:eigendec}
	\Kc = \KKI\otimes\KKO + \sig^2\III\otimes\IIO.
\end{align}
for moderately sized matrices $\KKI$ and $\KKO$ we  discuss how to use the full eigendecomposition
$\KKI = \UUI\LLI\UUI^\top$ and $\KKO = \UUO\LLO\UUO^\top$.

If we assume invertibility of $\KKO,$ we can write
\begin{align}
	\Kc &= \left((\KKI\otimes\IIO) + \sig^2(\III\otimes\KKO^{-1})\right)(\III\otimes\KKO).
\end{align}
Next, we insert the eigendecomposition and factorize
\begin{align}
	\Kc &= (\UUI\otimes\UUO)\left(\LLI\otimes\IIO + \sig^2\III\otimes\LLO^{-1}\right)(\UUI^\top\otimes\UUO^\top)(\III\otimes\UUO\LLO\UUO^\top).
\end{align}
By applying $(\CC_1\otimes\CC_2)^{-1}=\CC_1^{-1}\otimes\CC_2^{-1}$ and substitute $\HH=\LLI\otimes\IIO+\sig^2\III\otimes\LLO^{-1}$ we can write the inverse of $\Kc$ as 
\begin{align}
	\Kc^{-1}
	&= (\III\otimes\UUO\LLO^{-1}\UUO^\top) (\UUI\otimes\UUO) \HH^{-1} (\UUI^\top\otimes\UUO^\top).
\end{align}
The computation of the posterior mean $\myvec{\MMp}=\MUp$ becomes
\begin{align}
	\myvec{\MMp} &= \left({\KKI^\ast}^\top\otimes{\KKO^\ast}^\top\right)(\III\otimes\UUO\LLO^{-1}\UUO^\top) (\UUI\otimes\UUO) \HH^{-1} \myvec{\UUO^\top\Yu\UUI},
\end{align}
where $\yuv=\myvec{\Yu}$.
Further, we use the substitution $\myvec{\QQ}=\HH^{-1}\myvec{\UUO^\top\Yu\UUI},$ which has the equivalent expression
\begin{align}
	\QQ\LLI+\sig^2\LLO^{-1}\QQ=\UUO^\top\Yu\UUI.
\end{align}
Since $\LLI$ and $\LLO$ are diagonal matrices, this equation can be decoupled into independent scalar equations
\begin{align}
	(\QQ)_{ij} = \left( (\LLI)_j + \sigma^2(\LLO^{-1})_i \right)^{-1}(\UUO^\top\Yu\UUI)_{ij},
\end{align}
or by using the Hadamard product and Hadamard inverse
\begin{align}
	\QQ = \left(\mathbf{1}_{\Ic}\LLI + \sigma^2\LLO^{-1}\mathbf{1}_{\Oc} \right)^{\circ-1}\circ\left(\UUO^\top\Yu\UUI\right).
\end{align}
In conclusion we obtain
\begin{align}
	\MMp &= {\KKO^\ast}^\top\UUO\LLO^{-1}\QQ\UUI^\top{\KKI^\ast}^\top.
\end{align}
For large matrices $\KKI$ and $\KKO$  computing their  eigendecomposition is not feasible anymore and iterative solvers have to be used.

\section{Iterative solvers using Krylov subspaces}\label{sec:iterative_solvers}
In order to avoid the costly computation of the eigendecomposition of the possibly large matrices $\KKO$ and $\KKI$ we here derive a linear system solver not based on eigenfactorizations of the two covariance matrices. To compute the posterior mean of multiple outputs with covariance matrices $\KKI$, $\KKO$ and noise $\sig^2>0$ we saw in the previous sections that we have to solve
\begin{align}
	\left(\KKI\otimes\KKO+\sig^2\III\otimes\IIO\right)\myvec{\Xb}=\myvec{\Yu}
\end{align}
for $\Xb$ of fitting size. By using \Cref{eqn:kronecker_vec} we see, that this is equivalent to solving
\begin{align}\label{eqn:KKOstein}
	\KKO\Xb\KKI+\sig^2\Xb=\Yu.
\end{align}

This is a linear matrix equation \cite{bennerNumericalSolutionLarge2013,simonciniComputationalMethodsLinear2016} and if we divide \Cref{eqn:KKOstein} by $\sig^2$ we receive the Stein equation \cite{zhouSmithtypeIterativeAlgorithms2009}, which for general matrices $\AA$ and $\BB$ is denoted by
\begin{align}\label{eqn:stein}
	\AA\Xb\BB+\Xb=\CC.
\end{align} 

The Smith iteration method \cite{smithMatrixEquationXA1968} formally solves \Cref{eqn:stein} and is given by
\begin{align}
	\Xb_{k+1}=\AA\Xb_k\BB+\CC, \quad \Xb_0=\CC.
\end{align}
This method typically converges too slowly but can be accelerated in the following way
\begin{equation}\label{eqn:smith_accelleration}
	\begin{alignedat}{2}
		&\Xb_{k+1} = \AA_k\Xb_k\BB_k + \Xb_k, \\
		&\AA_k = \AA_{k-1}^2, \quad \BB_k = \BB_{k-1}^2, \\
		&\Xb_0 = \CC, \quad \AA_0 = \AA, \quad \BB_0 = \BB,
	\end{alignedat}
\end{equation}
see \cite{bennerNumericalSolutionLargescale2011, bennerSquaredSmithMethod2014, liLargescaleSteinLyapunov2013, ramadanHessenbergMethodNumerical2010a, sadkaneLowrankKrylovSquared2012, wuComplexConjugateMatrix2017, zhouSmithtypeIterativeAlgorithms2009}. Since we have large sparse matrices, explicitly computing their matrix powers like in \Cref{eqn:smith_accelleration} would lead to dense matrices and is therefore not feasible.

Assuming the invertibility of one of the kernel matrices, the system can be written as a Sylvester equation
\begin{align}\label{eqn:sylvester}
	\AA\Xb+\Xb\bar{\BB}=\bar{\CC},
\end{align}
for which a large body of literature exists; see, e.g., \cite{bennerNumericalSolutionLarge2013, simonciniNumericalSolutionOfAX1996} and the references therein. Note that the Stein equation and the Sylvester equation are tightly connected if we assume invertibility of at least one of the matrices $\AA$ or $\BB$. Note that in general, $\AA$ and $\bar{\BB}$ may be of different dimensions, resulting in a rectangular solution matrix $\Xb$. A unique solution for $\CC \ne 0$ is guaranteed if the spectra of $\AA$ and $-\bar{\BB}$ are disjoint. Importantly, when adapting Sylvester equation solvers to solve Stein equations, one should avoid explicitly computing matrix inverses and instead (approximately) solve with the matrix within the method. For small- to medium-scale problems (with matrix dimensions up to a few hundred), robust numerical methods such as the Bartels-Stewart algorithm \cite{bartelsAlgorithm432C21972}, which relies on the full Schur decomposition of the involved matrices, or the Hessenberg-Schur method \cite{golubHessenbergSchurMethodProblem1979}, can be effectively employed.

Since we mainly consider large $\KKI$ and $\KKO$, these methods are not feasible for our setup. To solve the Sylvester equation for large $\AA$ and $\bar{\BB}$, projections methods or the ADI method are the tools of choice. While we explain the projection method in this section for the \kpik~method in detail, we only briefly look at the ADI method.

The original derivation of the method was done with using the connection between Sylvester and Stein equations. Namely, we have that \Cref{eqn:sylvester} can be written as
\begin{align}
	(q\II+\AA)\Xb(q\II+\bar{\BB}) - (q\II-\AA)\Xb(\bar{\BB}-q\II)=2q\bar{\CC}, \quad q\neq 0.
\end{align}
Assuming that $q>0$ we can premultiply by $(q\II+\AA)^{-1}$ and postmultiply by $(q\II+\bar{\BB})^{-1}$ leads to
\begin{align*}
\Xb -\underbrace{(q\II+\AA)^{-1}(q\II-\AA)}_{\Ac}\Xb\underbrace{(\bar{\BB}-q\II)(q\II+\bar{\BB})^{-1}}_{\Bc}=\underbrace{2q(q\II+\AA)^{-1}\CC(q\II+\bar{\BB})^{-1}}_{{\Cc}}.
\end{align*}
We can now efficiently write this as another Stein equation (cf. \cite{simonciniComputationalMethodsLinear2016,smithMatrixEquationXA1968})
\begin{align}
	\Xb-\Ac\Xb\Bc = \Cc
\end{align}
for which the solution can formally be expressed as a series expansion, which under the condition that the spectral radius of the matrices $\Ac$ and $\Bc$ is less than one converges. If one introduces two parameters, e.g. $p$ and $q$ one obtains another two-parameter Stein equation. Again the condition of both matrices $\Ac$ and $\Bc$ having a spectral radius less than one guarantees the existence of a solution. This has led to hunt for the optimal \textit{shift parameters} $p$ and $q$ with minimal spectral radii such that  fast convergence of the numerical methods is obtained (cf. \cite{bennerADIMethodSylvester2009a, simonciniComputationalMethodsLinear2016}). The review articles \cite{bennerNumericalSolutionLarge2013,simonciniComputationalMethodsLinear2016} provide excellent starting points for further details and more recent developments.

We here focus on Krylov subspace based methods. Namely, we derive a low-rank preconditioned \cg\ (\lrpcg) method \cite{hestenesMethodsConjugateGradients1952}. For this we use the Krylov-plus-Inverse-Krylov method (\kpik) developed in \cite{simonciniNewIterativeMethod2007} as preconditioner, which is another Krylov subspace method.

\subsection*{The \lrpcg\ solver}
We briefly recall the main steps of the \cg\ method and illustrate how these can be performed in low-rank fashion. We will later comment on the use of the matrix equation solvers for the corresponding Sylvester equation. In \Cref{alg:pcg} the principal preconditioned \cg\ is shown.

\begin{algorithm}
	\caption{Preconditioned Conjugate Gradient (PCG) Method}
	\label{alg:pcg}
	\begin{algorithmic}[1]
		\REQUIRE spd matrix $\Acg$, and right hand side $\bcg,$ possibly preconditioner $\Mcg$
		\STATE $\rr_0 = \bcg - \Acg\xx_0$
		\STATE Solve $\Mcg\zz_0 = \rr_0$
		\STATE $\pp_0 = \zz_0$
		\STATE $k = 0$
		\WHILE{not converged}
		\STATE $\alpha_k = \rr_k^\top \zz_k/(\pp_k^\top \Acg\pp_k)$
		\STATE $\xx_{k+1} = \xx_k + \alpha_k\pp_k$
		\STATE $\rr_{k+1} = \rr_k - \alpha_k\Acg\pp_k$
		\STATE  Solve $\Mcg\zz_{k+1} = \rr_{k+1}$
		\STATE $\beta_k = \zz_{k+1}^\top \rr_{k+1}/(\zz_k^\top \rr_k)$
		\STATE $\pp_{k+1} = \zz_{k+1} + \beta_k\pp_k$
		\STATE $k = k + 1$
		\ENDWHILE
	\end{algorithmic}
\end{algorithm}

In order to use this method for low-rank computation we need to make sure that every step can be performed maintaining this structure. We briefly illustrate this on the matrix vector product, vector updates, and the computation of an inner product. We focus on 
\begin{align}
	\pp=\myvec{\UU_{\pp}\VV_{\pp}^\top}
\end{align}
where we drop the iteration index $k$. Then the matrix vector product $\Acg\pp$, with $\Acg=\Kc$ from \cref{eqn:K} becomes in factorized form
\begin{align}
	\KKI \UU_{\pp}\VV_{\pp}^\top\KKO + \sig^2\UU_{\pp}\VV_{\pp}^\top=
	\begin{bmatrix}
		\KKI \UU_{\pp}&	\sig\UU_{\pp}\\
	\end{bmatrix}
	\begin{bmatrix}
		\KKO\VV_{\pp}&
		\sig\VV_{\pp}\\
	\end{bmatrix}^\top
\end{align}
where we have used the symmetry of $\KKO$. This will be used as $\Amult$ in the \lrpcg\ algorithm. We see that now the result is still a low-rank matrix of rank $2k$ if $\UU_{\pp}\VV_{\pp}^\top$ was a rank $k$ approximation. This rank growth can usually be truncated by following every potential rank increase by a truncation procedure, usually based on an economic QR decomposition or a truncated SVD \cite{kressnerLowRankTensorKrylov2011,stollLowRankTimeApproach2015}. The update of a vector, e.g., $\xx_{k+1} = \xx_k + \alpha_k\pp_k$ represented in low-rank form as
\begin{align}
	\begin{bmatrix}
		\UU_{\xx}&	\alpha \UU_{\pp}\\
	\end{bmatrix}
	\begin{bmatrix}
		\VV_{\xx}&\VV_{\pp}\\
	\end{bmatrix}^\top
\end{align}
again with some rank growth that is controlled by a truncation step. This will be called $\lrsum$ in the \lrpcg\ algorithm. For the computation of inner products and norms we rely on the fact that it holds for two matrices $\SS_1, \SS_2$
\begin{align}
	\trace{\SS_1^\top \SS_2}=\myvec{\SS_1}^\top\myvec{\SS_2}
\end{align}
as using this leads to
\begin{align}
	\trace{(\UU_{\rr}\VV_{\rr}^\top)^\top(\UU_{\zz}\VV_{\zz}^\top)}=
	\trace{(\UU_{\rr}^\top\UU_{\zz})(\VV_{\zz}^\top\VV_{\rr})},
\end{align}
where ensuring that the ranks of the low-rank approximations match. This means the trace has to be computed for the product of two small matrices. We call this $\traceprod$ in the \lrpcg\ method. We then have the complete \lrpcg\ version in \cref{alg:lrpcg}.
\begin{algorithm}
	\caption{Low-Rank Preconditioned Conjugate Gradient (\lrpcg) Method. The index $k$ is dropped for all iterates. Updated quantities are indicated by a bar over the variable.}
	\label{alg:lrpcg}
	\begin{algorithmic}[1]
		\REQUIRE spd matrices $\AA_1, \AA_2$, and low-rank right hand side $\UU_{\BB}, \VV_{\BB}$, possibly preconditioner $\precon$
		\STATE $\UU_{\QQ}, \VV_{\QQ} = \Amult(\UU_{\Xb}, \VV_{\Xb})$
		\STATE $\UU_{\RR}, \VV_{\RR} = \lrsum(\UU_{\BB}, \VV_{\BB}, -\UU_{\QQ}, \VV_{\QQ})$
		\STATE $\UU_{\ZZ}, \VV_{\ZZ} = \precon(\bar \UU_{\RR}, \bar \VV_{\RR})$
		\STATE $\UU_{\PP}, \VV_{\PP} = \UU_{\ZZ}, \VV_{\ZZ}$
		\WHILE{not converged}
		\STATE \label{lst:line:amult} $\UU_{\TT}, \VV_{\TT} =\Amult(\UU_{\PP}, \VV_{\PP})$
		\STATE $\alpha = \traceprod(\UU_{\RR}, \VV_{\RR}, \UU_{\ZZ}, \VV_{\ZZ})/\traceprod(\UU_{\PP}, \VV_{\PP}, \UU_{\TT}, \VV_{\TT})$
		\STATE $\bar \UU_{\Xb}, \bar \VV_{\Xb} = \lrsum(\UU_{\Xb}, \VV_{\Xb}, \UU_{\PP}, \alpha \VV_{\PP})$
		\STATE $\bar \UU_{\RR}, \bar \VV_{\RR} = \lrsum(\UU_{\RR}, \VV_{\RR}, -\alpha \UU_{\TT}, \VV_{\TT})$
		\STATE $\bar \UU_{\ZZ}, \bar \VV_{\ZZ} = \precon(\bar \UU_{\RR}, \bar \VV_{\RR})$
		\STATE $\beta = \traceprod(\bar \UU_{\RR},\bar \VV_{\RR},\bar  \UU_{\ZZ}, \bar \VV_{\ZZ})/\traceprod(\UU_{\RR}, \VV_{\RR}, \UU_{\ZZ}, \VV_{\ZZ})$
		\STATE $\bar \UU_{\PP}, \bar \VV_{\PP} = \lrsum(\bar \UU_{\ZZ},\bar  \VV_{\ZZ}, \beta \UU_{\PP}, \VV_{\PP})$
		\ENDWHILE
		\STATE $\Amult(\UU, \VV) := \trunc([\sqrt{d}\UU, \AA_2\UU] , [\sqrt{d}\VV, \AA_1\VV])$
		\STATE $\lrsum(\UU_1, \VV_1, \UU_2, \VV_2) := \trunc\left(\begin{bmatrix} \UU_1 & \UU_2 \end{bmatrix}, \begin{bmatrix}\VV_1 & \VV_2 \end{bmatrix}\right)$
		\STATE $\traceprod(\UU_1, \VV_1, \UU_2, \VV_2) := \trace{(\UU_1^\top\UU_2) (\VV_2^\top\VV_r))}$
	\end{algorithmic}
\end{algorithm}

While Krylov methods are a method of choice in many application their convergence behaviour depends on the matrix properties such as the distribution of the eigenvalues \cite{liesenKrylovSubspaceMethods2013}. For a better performance of the solvers preconditioning is essential and we now discuss the choice of the preconditioner in the \lrpcg\ method. We introduce \kpik\ in some detail in the next section, which can serve as a preconditioner but more prominently as a standalone solver. 

\subsection*{The \kpik\ solver}
Here, we already emphasized that we consider the case where the coefficient matrices $\AA$ and $\BB$ are both large, possibly dense, and computing the spectral decompositions and backward solves is too computationally expensive. As already assumed before we want to assume that $\CC$ has much smaller rank than the problem dimensions of $\AA$ or $\BB$. This results in low-rank methods for Sylvester equations where the \lrpcg\ method from the previous section is one of those examples. These low-rank approximations are of special interest, since in general $\Xb$ is dense, and thus hard to store when $\AA$ and $\BB$ are sufficiently large.

For \kpik\ we consider the projection method based on two approximation spaces $\textrm{Range}(\VV)$ and $\textrm{Range}(\WW)$, an approximation $\Xbt = \VV \YYt \WW^T$ is determined by requiring that the residual matrix $\RR:= \AA \Xbt + \Xbt \BB - \CC$ satisfies the Galerkin orthogonality conditions
\begin{align}
	\VV^T \RR \WW = 0.
\end{align}
If both $\VV$ and $\WW$ have orthogonal columns and using the fact that $\Xbt = \VV \YYt \WW^T$ we get the reduced matrix equation
\begin{align}
	(\VV^T \AA \VV) \YYt + \YYt (\WW^T \BB \WW) - \VV^T \CC \WW = 0,
\end{align}
where now $\VV^T \AA \VV$ and $\WW^T \BB \WW$ are of small dimensionality. The reduced equation can now be solved efficiently using the Bartels-Stewart method providing us with $\YYt$. The choice of Range($\VV$) and Range($\WW$) leads to different approximate solutions (cf. \cite{simonciniComputationalMethodsLinear2016})

Our choice here is for Range($\VV$) to correspond to an approximation space generated by powers of $\AA$ and $\AA^{-1}$ and it was first proposed under the name of {\it Extended Krylov subspace} \cite{druskinExtendedKrylovSubspaces1998}. In \cite{simonciniNewIterativeMethod2007} it was shown for $\BB=\AA^T$ that such a space can be generated as
\begin{align}
	{\mathbb E}{\mathbb K}(\AA, \CC_1)=
	{\rm Range}([\CC_1, \AA^{-1}\CC_1, \AA \CC_1, \AA^{-2}\CC_1, \AA^2 \CC_1, \AA^{-3}\CC_1, \ldots ])
\end{align}
and expanded until the approximate solution $\Xbt$ is sufficiently good.  Note that in a standard implementation that sequentially generates ${\mathbb E}{\mathbb K}(\AA, \CC_1)$, two blocks of vectors are added at each iteration, one for $\AA$, and one for $\AA^{-1}$.  The implementation of the resulting projection method with ${\mathbb E}{\mathbb K}(\AA, \CC_1)$ as approximation space was called \kpik\ in \cite{simonciniNewIterativeMethod2007} for the Lyapunov matrix equation, that is $\BB=\AA^T$ and $\CC_1=\CC_2$. The procedure in \cite{simonciniNewIterativeMethod2007} was adapted to the case of the Sylvester equation with inexact solves in \cite{breitenLowRankSolversFractional2016}, via the space ${\mathbb E}{\mathbb K}(\BB^T, \CC_2)$, which we will refer to as \kpik\ or once used with inexact solves \ikpik.	

From a computational standpoint, our application problem is particularly demanding because the	iterative generation of the extended space requires solving systems with $\AA$ and $\BB$, whose size can be very large. Furthermore we recall that we are interested in the Stein \Cref{eqn:stein1} but that this can be transformed to a Sylvester equation if the kernel matrices are invertible. We then obtain the system 
\begin{equation}
	\label{eqn:sylvstein}
	\Xb\KKO + \sig^2 \KKI^{-1}\Xb=\Yu.	
\end{equation}
This means that in the Sylvester framework we have $\AA=\sig^2 \KKI^{-1}$ and $\BB=\KKO$. Since we are using \kpik\ as a solver we need to solve \textit{and} multiply with the matrices $\KKO$ and $\KKI$ anyways, so this transformation does not increase the computational cost of \ikpik . We can now use the \ikpik\ as a standalone solver but also as a preconditioner in \cref{alg:lrpcg}.

\section{Numerical experiments}\label{sec:numerical_experiments}
Our aim is to compare the performance of the presented methods, namely the \lrpcg- and \ikpik~method. For our experiments, we adopted a model variant from \cite{zhiGaussianProcessesGraphs2023}.
\begin{figure}
	\centering
	\caption{Oxford street network example: 30\% of randomly chosen input nodes and 70\% output nodes (red).}
	\label{fig:ioOxford}
	\includegraphics[height=8cm]{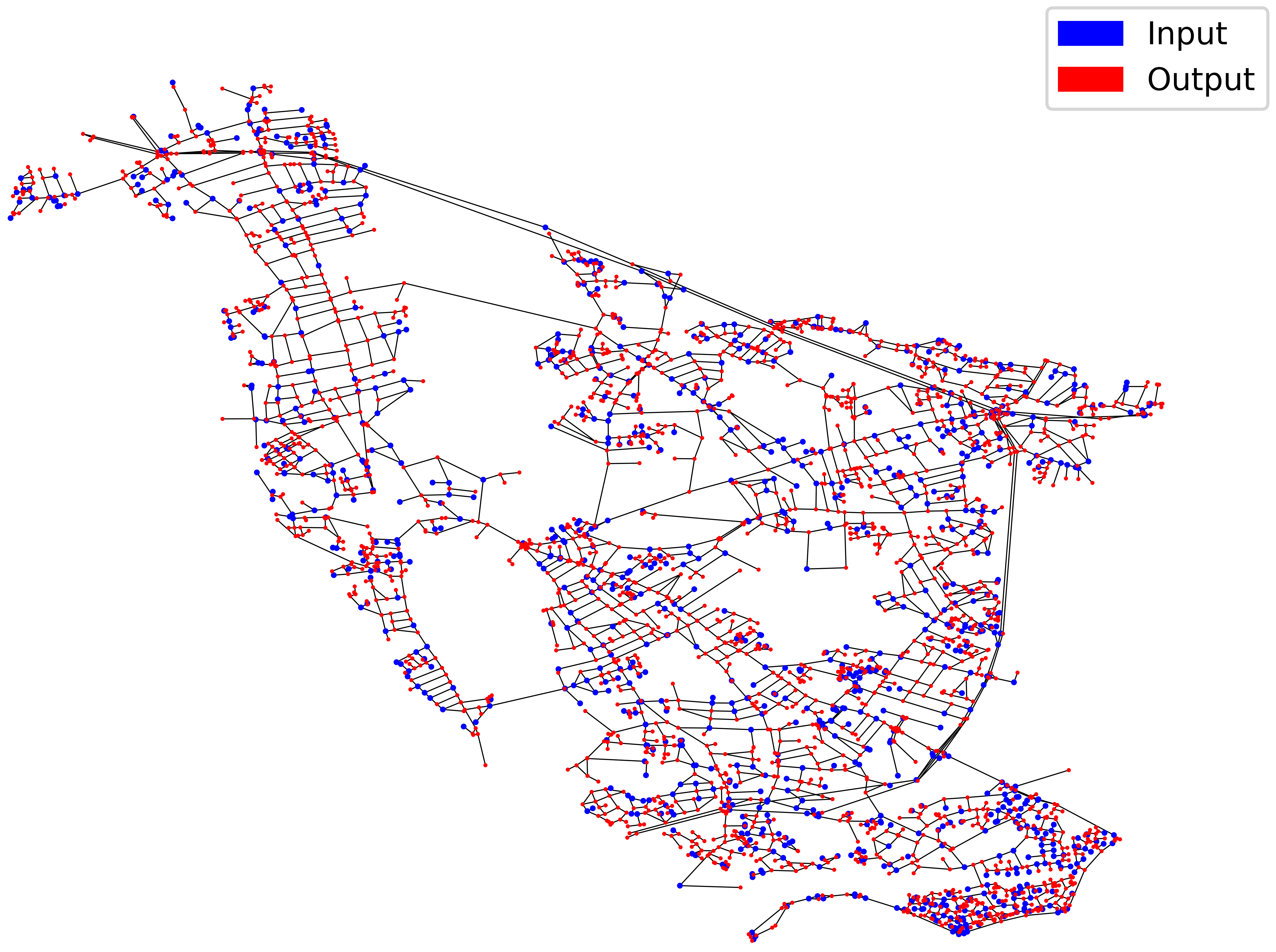}
\end{figure}

Assume that we have data, i.e. signals, for each time point in $\Tc$ given on the nodes of a graph $\Gc=(\Wc, \Ec)$. For another set of time points $\Tc^\ast$ we just know the signals on a subset $\Uc\subset\Wc$ but we want to infer them on the remaining nodes $\Uc^\ast=\Wc\setminus\Uc$. The idea to complete the missing information is: if values are similar on $\Uc$ at a given point in time, the corresponding values on $\Uc^\ast$ should also be similar. Therefore, we define the covariance between signals on $\Uc^\ast$ by a covariance function that is dependent on the distance of signals on $\Uc$. This way, similar input signals--even if far apart in time--lead to similar output signals, which makes the covariance independent of time. We call $\Uc$ the input nodes and $\Uc^\ast$ the output nodes with an illustrated example in \Cref{fig:ioOxford}. In our experiments we split them in ratio 1 to 4, e.g. 20\% are input nodes. Since we assume here, that the output training and target nodes are equal, i.e. ($\Vc=\Vc^\ast=\Uc^\ast$), it holds $\KKO=\KKO^\ast=\KKO^{\ast\ast}$. In \Cref{sec:submatrices} we give an idea for implementing the general case for $\Vc\neq\Vc^\ast$.

For our experiments, we use real-world street network graphs from the \texttt{osmnx} package, which where adjusted to be undirected, loop-free, and connected. On such a graph, we generate synthetic data using a dynamical system, i.e. an Allen-Cahn equation, as described in more detail next.

\subsection*{Synthetic data using Allen-Cahn-equation}
For the synthetical data $\DD$, we use a dynamical system which evolves on the graph over the time points, i.e. the Allen-Cahn-equation \cite{bartelsAllenCahnEquation2015}
\begin{align}
	\frac{\partial u}{\partial t} = \eps D \Delta u + \frac{1}{\eps}(u-u^3), 
\end{align}
with interface parameter $\eps\in\RE$ and diffusion parameter $D\in\RE$. In discrete form on the graph for time step size $\tau$ it becomes
\begin{align}\label{eqn:allencahndiscrete}
	(\II+\tau\eps D \LL) \uu_{k+1} = \uu_k + \frac{\tau}{\eps}(\uu_k-\uu_k^3).
\end{align}

Here, $\uu_k^3$ denotes the element wise power. At first, we perform LU decomposition of $\II+\tau\eps D \LL$ and compute step by step $\uu_2$ to $\uu_n$ by solving \Cref{eqn:allencahndiscrete} with parameters $\eps=0.08$, $D=100$, $\tau=5\cdot 10^{-4}$. We randomly initialize $\uu_1$.

The resulting data matrix $\DD$ is row-indexed in order $[\Uc, \Uc^\ast]$ and column-indexed in order $[\Tc, \Tc^\ast]$, giving the block structure
\begin{align}
	\DD = \begin{bmatrix}
		\Xu & \Xu^\ast \\
		\Yu & \Yu^\ast
	\end{bmatrix}.
\end{align}
In accordance with \Cref{eqn:mogp_matrices}, $\Xu, \Xu^\ast$ are the inputs to corresponding outputs $\Yu, \Yu^\ast$, where we assume low-rank factorization $\Yu = \UU_{\Yu} \VV_{\Yu}^\top$. The ground truth $\Yu^\ast$ is considered unknown. Since performing the regression only at $\Tc^\ast$  would lead to gaps in the regression plot, we chose the full time domain $\Tc\cup\Tc^\ast$ as targets in our experiment. Because we also want to use the full graph information, we take $\Uc\cup\Uc^\ast$ for the target. Hence, we extract it in the following way from the synthetic data
\begin{equation}
	\begin{aligned}
		\tilde{\Xu}       &= \DD_{\Uc, \Tc}, 
		&\quad \tilde{\Yu}       &= \DD_{\Uc \cup \Uc^\ast, \Tc}, \\
		\tilde{\Xu}^\ast  &= \DD_{\Uc, \Tc \cup \Tc^\ast}, 
		&\quad \tilde{\Yu}^\ast  &= \DD_{\Uc \cup \Uc^\ast, \Tc \cup \Tc^\ast}.
	\end{aligned}
\end{equation}
The matrix $\DD_{\Uc, \Tc}$ corresponds to selecting the rows $\Uc$ and columns $\Tc$ from $\DD$ and analogously for the other matrices.

\subsection*{Covariance matrices}
We test the solvers for large matrices and for different conditions of the system
\begin{align}\label{eqn:cov_system}
	\KKI \otimes \KKO + \sigma^2 \III\otimes\IIO.
\end{align}
For the input covariance matrices, we use a squared exponential kernel \cref{eqn:se_kernel} with $\sig^2_w=1$ over a time interval of total length 10,000, which we split by a ratio of 1 to 9 into $\Tc$ and $\Tc^\ast$, resulting, in an $1,000 \times 1,000$ matrix $\KKI$. For the output covariance matrix $\KKO$, we implemented both the global filter $\BB=(\II+\alpha \LL)^{-1}$ and the local average filter $\BB=(\II+\alpha \DD)^{-1}(\II+\alpha \WW)$. We construct these based on the adjacency matrices of the real-world street network graph of Berlin with 28,189 nodes and 42,265 edges. Since the filter matrices are defined as the inverses of sparse matrices, we have to find efficient strategies to multiply and solve with them. To ensure efficiency, we employ the sparse linear algebra libraries \texttt{scipy} (for linear operator arithmetic) and \texttt{sksparse} (for sparse Cholesky decomposition).

The global filter matrix is given by
\begin{align}
	\KKO = (\II+\alpha \LL)^{-2}.
\end{align}
Multiplication by $\KKO$ is solving by $(\II+\alpha \LL)^2$. Since $\LL$ has non-negative eigenvalues, for $\alpha>0$ the matrix $\II+\alpha \LL$ is positive definite. Hence, we suppose Cholesky decomposition of $\II+\alpha \LL$ and solve twice by the Cholesky factors. Solving with $\KKO$ becomes very efficient by multiply twice by $\II+\alpha \LL$.

For the local average filter we obtain
\begin{align}
	\KKO=(\II+\alpha \DD)^{-1}(\II+\alpha \WW)^2(\II+\alpha \DD)^{-1}.
\end{align}
Since $\II+\alpha \WW$ is in general not positive definite we explicitly form $\KKO$ and compute its Cholesky decomposition. The main cost comes from computing $(\II+\alpha \WW)^2$ as the other terms are diagonal.

To control the condition of \Cref{eqn:cov_system}, we vary three hyperparameters: the length scale $\ell > 0$, which affects the condition number of $\KKI$; the regularization parameter $\alpha > 0$, which influences the condition number of $\KKO$; and the noise variance $\sigma^2$, which governs the overall conditioning of the system. If $\lambda_O$ and $\mu_O$ are the largest and smallest eigenvalues of $\KKO$ and $\lambda_I$ and $\mu_I$ the largest and smallest eigenvalue  of $\KKI$ respectively, the total condition of the system is $$\frac{\lambda_O\lambda_I+\sig^2}{\mu_O\mu_I+\sig^2}<1+\lambda_O\lambda_I/\sig^2.$$ For our experiments, we choose the default hyperparameters values as $\ell = 10$, $\alpha = 1$, and $\sig^2 = 5 \cdot 10^{-3}$, and vary only one at a time.

The framework constructs the covariance matrices according to the selected hyperparameters and then solves the Stein equation by the proposed methods. Both \lrpcg~and \ikpik~terminate if a relative residual below $\varepsilon = 10^{-8}$ is reached.
For \lrpcg, we apply a truncation tolerance of $10^{-10}$ and use 2 steps of \ikpik~for preconditioning.

\subsection*{Degree-Weighted Average}
Additionally, we want to test the iterative solvers for the covariance matrix discussed in \Cref{sec:degree_weighted_average}, which leads to an efficient way of solving \Cref{eqn:cov_system}. To demonstrate this we use the Germany street network graph with 3,025,346 nodes and 3,913,844 edges. As we assumed in \Cref{sec:degree_weighted_average}, this covariance matrix is suitable for data with a stationarity property, e.g. it is (approximately) invariant under the application of the transition matrix $\PP=\DD^{-1}\WW$. Therefore, we compute synthetical data following the precise stationary property by first computing the stationary distribution
\begin{align}
	\ss_{k+1}=\PP^\top\ss_k, \quad \ss_0=\ee_1\in\RE^w,
\end{align}
with $\ee_1$ the first unit vector. We stop the iteration when $\norm{\ss_{k+1}-\ss_k}<10^{-6}$, where we call the last iterate $\st$. From this, we construct the data like
\begin{align}
	\DD = \underbrace{\begin{bmatrix}
		\st & \dots & \st
	\end{bmatrix}}_{10,000 \text{ times}} + \EPS,
\end{align}
where we add a normally distributed noise term $\EPS$ with size $w\times 10,000,$ whose entries are all i.i.d. We would like to emphasize that our primary goal is to test the methods even though the data are only linear over time. However, as we have shown in \Cref{sec:degree_weighted_average}, the methods also work for complex data as long as it satisfies the stationary assumption.

We chose the squared exponential kernel for input covariance and the output covariance according to \Cref{eqn:stat_cov_matrix} where we don't have to decide for $\SIGMAp{22}$, since it is not necessary for computing the posterior mean. In contrast to the previous model for the filter matrices, we split the nodes into training and target nodes, i.e., $\Vc \neq \Vc^\ast$. To avoid gaps in the regression plot, we use again the complete time interval $\Tc\cup\Tc^\ast$ for the target. The data is used in the following way
\begin{equation}
	\begin{aligned}
		\tilde{\Xu}       &= \DD_{\Vc, \Tc}, 
		&\quad \tilde{\Yu}       &= \DD_{\Vc, \Tc}, \\
		\tilde{\Xu}^\ast  &= \DD_{\Vc, \Tc \cup \Tc^\ast}, 
		&\quad \tilde{\Yu}^\ast  &= \DD_{\Vc^\ast, \Tc \cup \Tc^\ast}.
	\end{aligned}
\end{equation}
Since the degree-weighted average is independent of hyperparameters we did not vary $\alpha$.

\subsection*{Hardware Specifications} All computations were performed on a machine equipped with 2 × AMD EPYC 9534 64-Core processors. The machine is equipped with approximately 1.5 TB of RAM.

\subsection{Results}

\begin{figure}
	\centering
	\begin{minipage}{0.48\textwidth}
		\centering
		\caption{Regression on data from dynamical system}
		\label{fig:regression}
		\includegraphics[height=4.4cm]{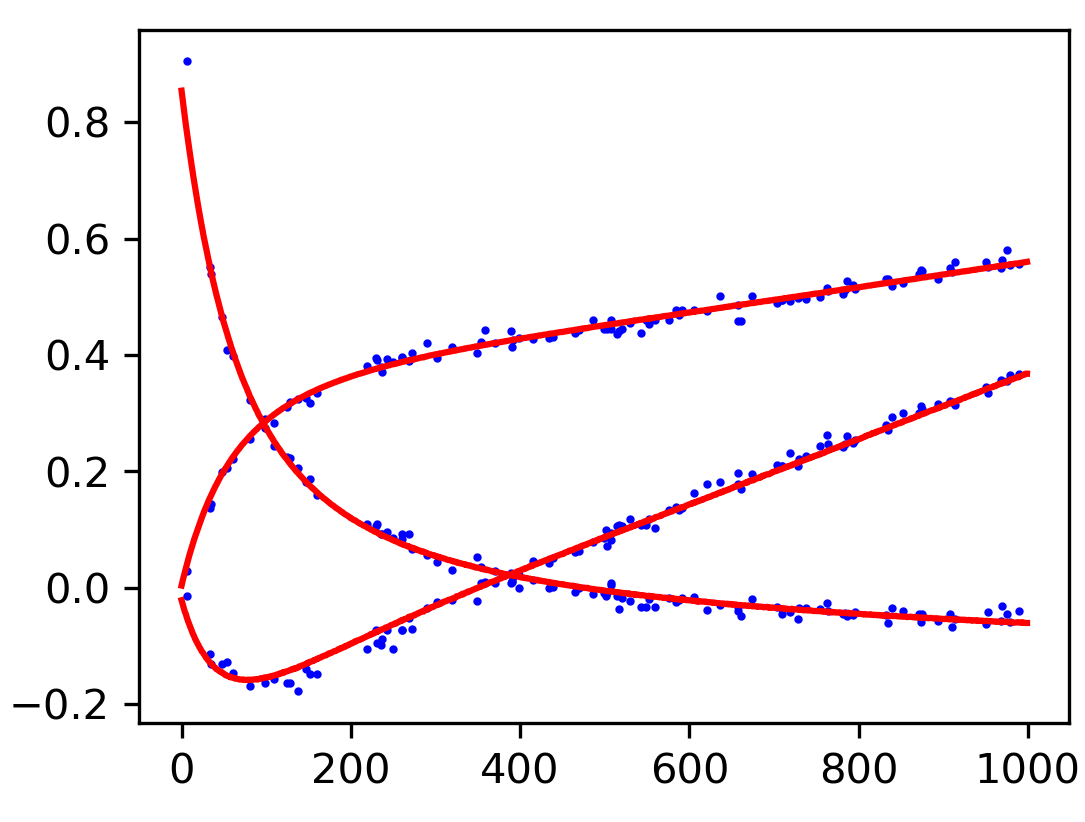}
	\end{minipage}
	\hfill
	\begin{minipage}{0.48\textwidth}
		\centering
		\caption{Regression on data following stationary property}
		\label{fig:regression2}
		\includegraphics[height=4.5cm]{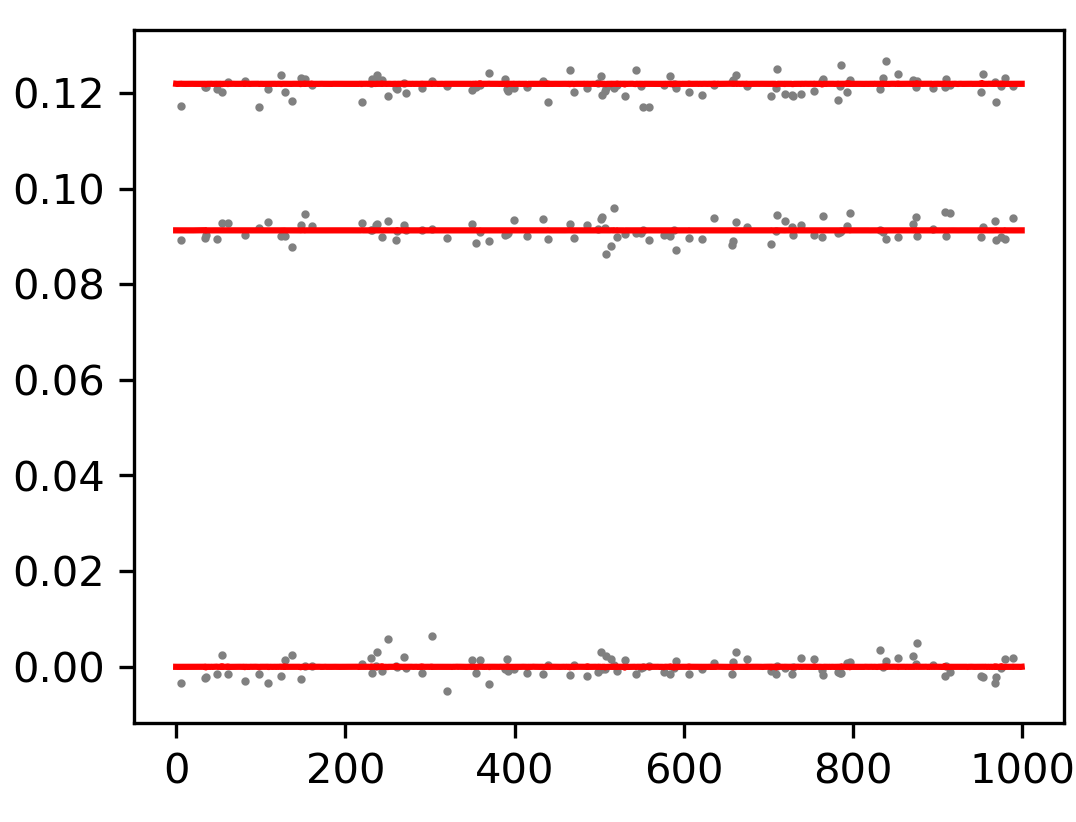}
	\end{minipage}
\end{figure}

\begin{figure}
	\begin{subfigure}[b]{0.47\textwidth}
		\caption{Iteration count for increasing $\alpha$}
		\label{fig:result1}
		\includegraphics[height=4.5cm]{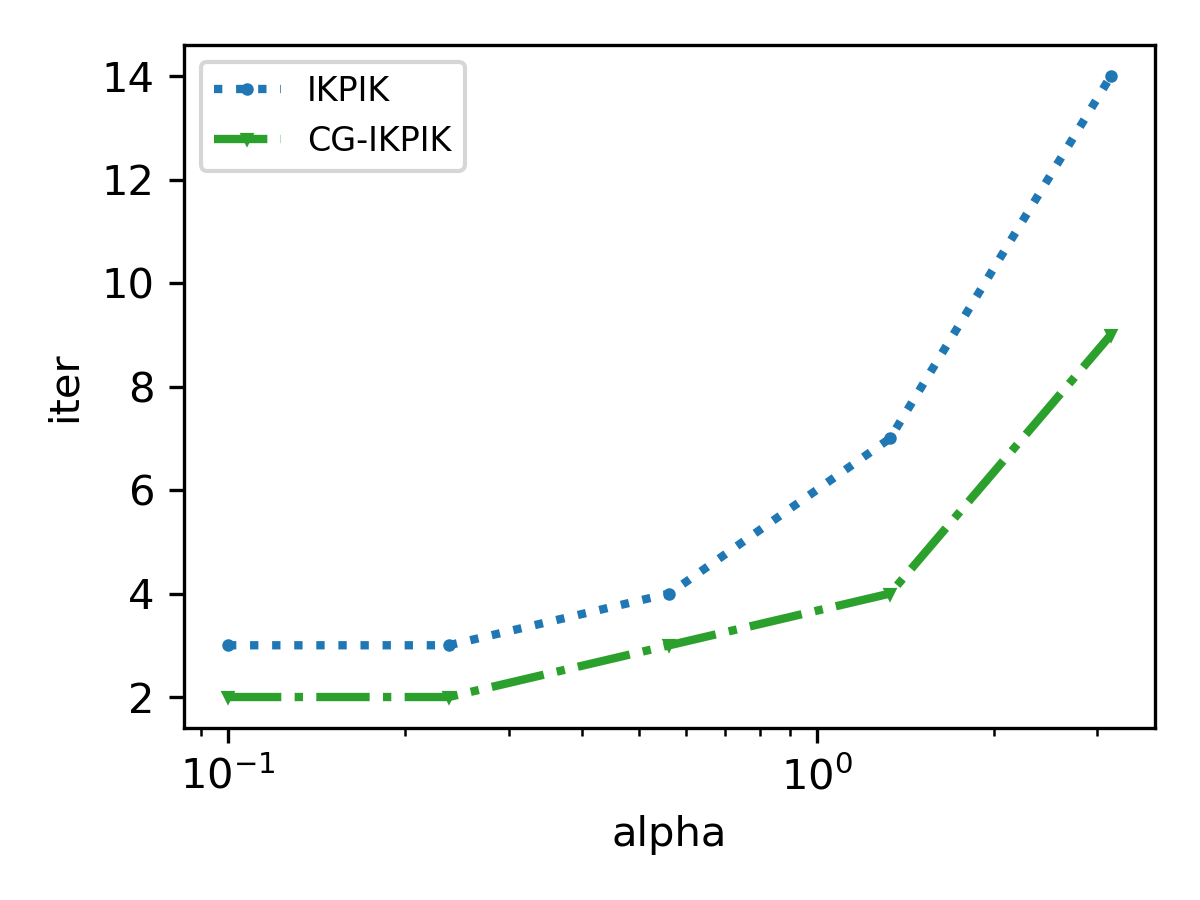}
	\end{subfigure}
	\begin{subfigure}[b]{0.47\textwidth}
		\caption{Sol. rank for increasing $\alpha$}
		\label{fig:result2}
		\includegraphics[height=4.5cm]{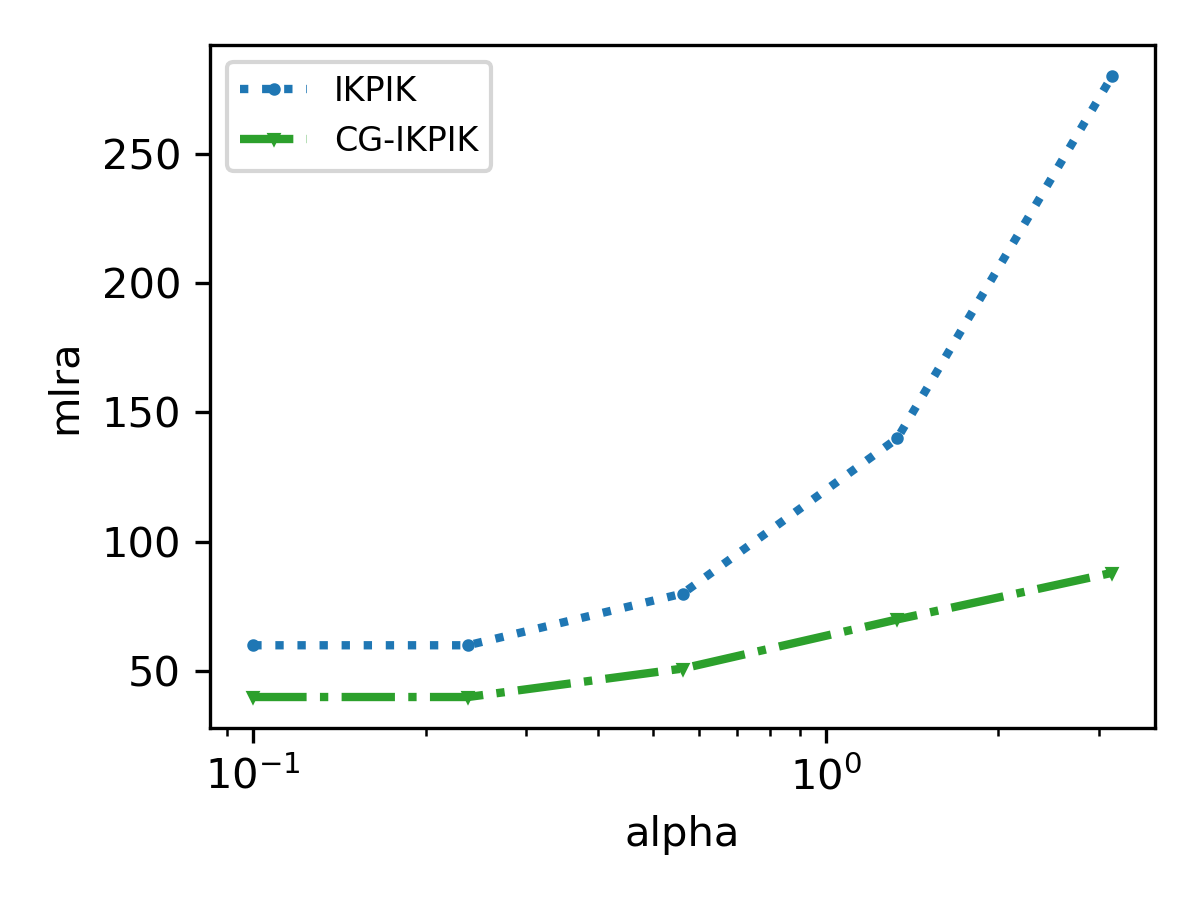}
	\end{subfigure}
	
	\begin{subfigure}[b]{0.47\textwidth}
		\caption{Runtime for increasing $\alpha$}
		\label{fig:result3}
		\includegraphics[height=4.5cm]{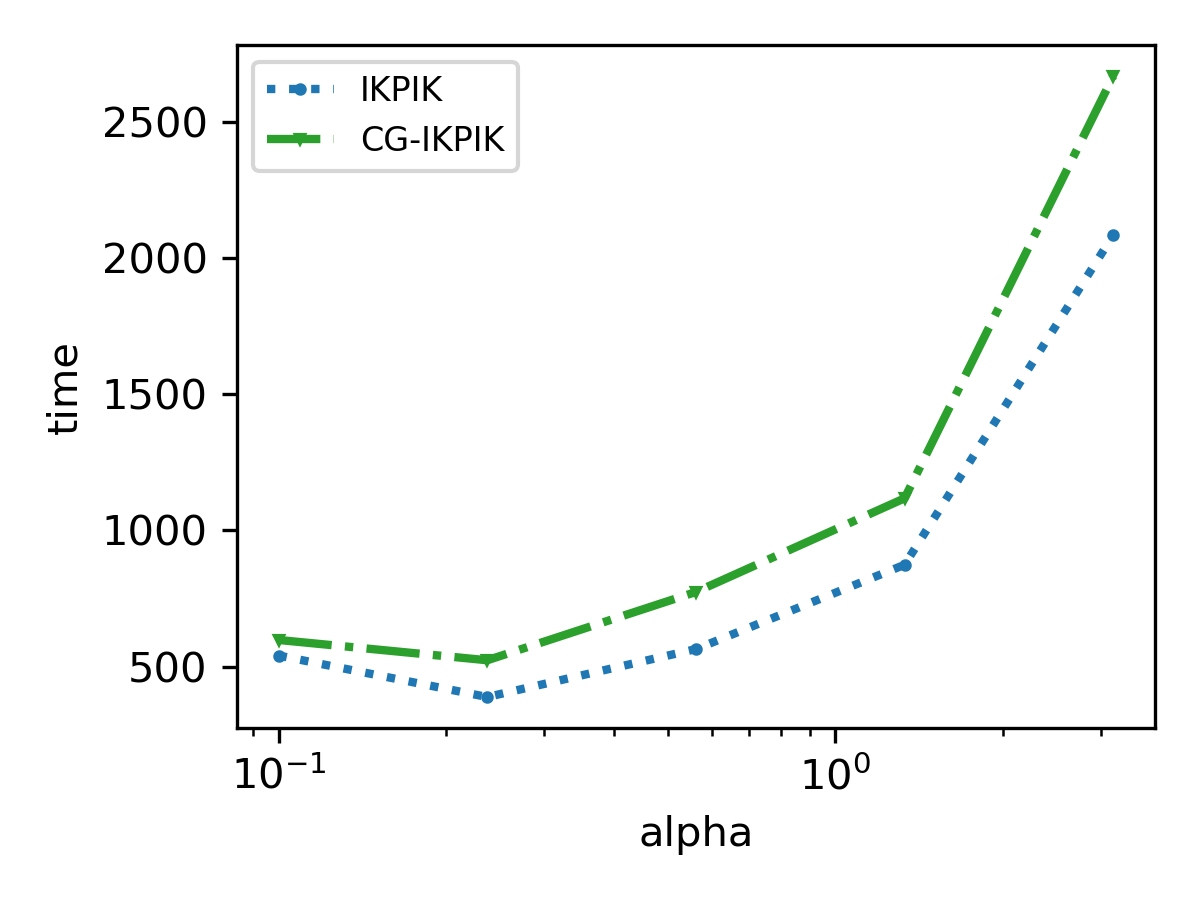}
	\end{subfigure}
	\begin{subfigure}[b]{0.47\textwidth}
		\caption{Residual norm for increasing $\alpha$}
		\label{fig:result4}
		\includegraphics[height=4.5cm]{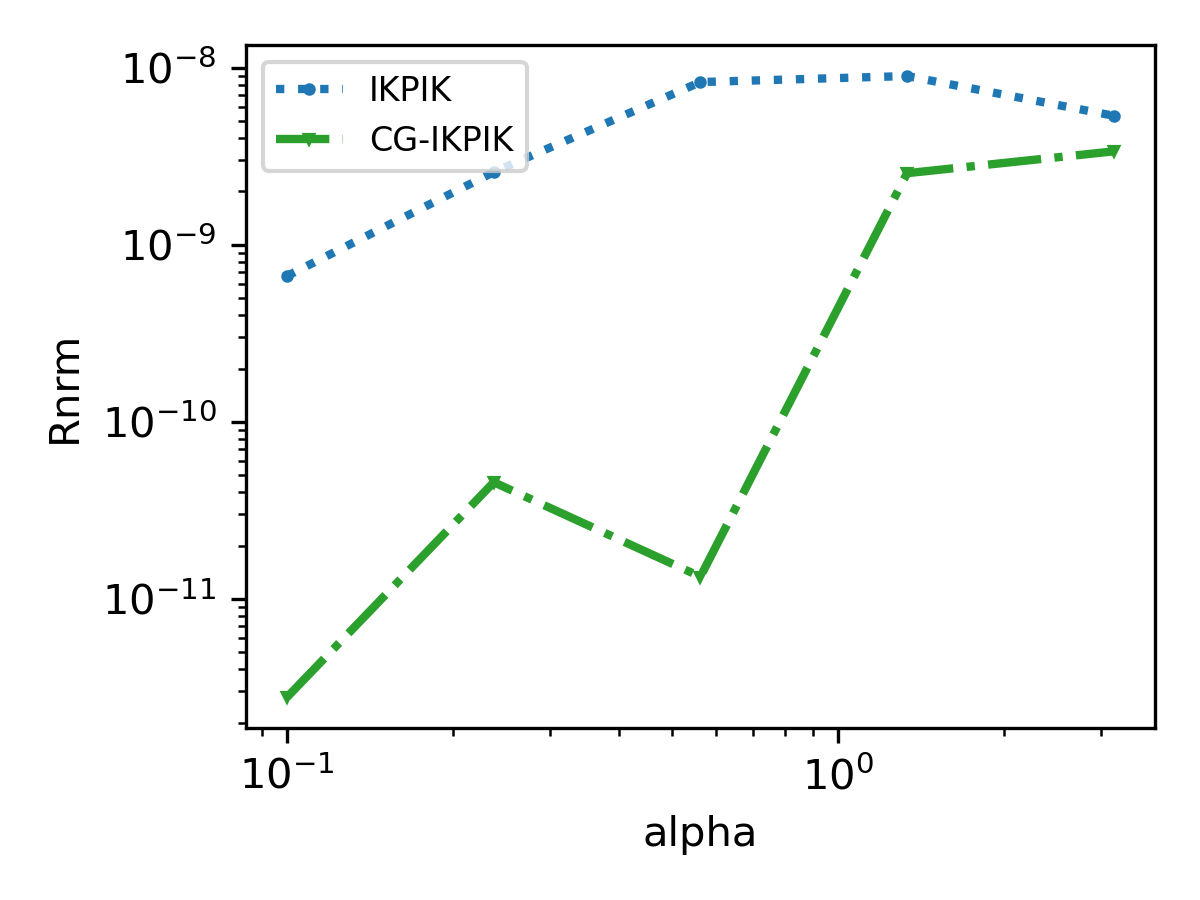}
	\end{subfigure}
	\caption{Comparison of the methods performance when varying $\alpha$.}
\end{figure}

\Cref{fig:regression,fig:regression2} show the regression plots for the global filter applied to data from a dynamical system, and for the degree-weighted average filter applied to data following the stationary distribution. As shown, the regression lines fit the data well. The gray data points are assumed to be unknown, but we are able to infer them based on information from neighboring nodes. This kind of information transfer would not be possible in a setup using only independent single GPs. In \Cref{fig:result1,fig:result2,fig:result3,fig:result4} iteration count, solution rank, runtime (in seconds) and residual norm of both methods for the Berlin street network with global filter are compared for increasing $\alpha$.

\begin{table}[h]
	\centering
	\caption{Global filter}
	\label{tab:globalfilter}
	\begin{tabular}{ c | r r  r r  r r}
		\toprule
		\textbf{global filter} & \multicolumn{2}{c}{\textbf{iter}} & \multicolumn{2}{c}{\textbf{sol. rank}} & \multicolumn{2}{c}{\textbf{runtime (s)}} \\
		& \ikpik & \lrpcg & \ikpik & \lrpcg & \ikpik & \lrpcg \\
		\midrule
		$\alpha=0.1$     &     2 &     1 &    40 &    31 &    362 &   318 \\
		$\alpha=2$		 &     3 &     2 &    60 &    33 &    447 &   545 \\
		$\alpha=10$		 &     5 &     3 &   100 &    23 &    650 &   780 \\
		\midrule
		$\ell=1$		 &     5 &     3 &   100 &    89 &    566 &   700 \\
		$\ell=5$		 &     3 &     2 &    60 &    49 &    334 &   475 \\
		$\ell=10$		 &     2 &     1 &    40 &    31 &    236 &   237 \\
		\midrule
		$\sig^2=10^{-2}$ &     3 &     2 &    60 &    32 &    400 &   480 \\
		$\sig^2=10^{-3}$ &     3 &     2 &    60 &    50 &    379 &   491 \\
		$\sig^2=10^{-4}$ &     4 &     3 &    80 &    79 &    469 &   726 \\
		\bottomrule
	\end{tabular}
\end{table}
\begin{table}[h]
	\centering
	\caption{Local filter}
	\label{tab:localfilter}
	\begin{tabular}{ c | r r  r r  r r}
		\toprule
		\textbf{local filter} & \multicolumn{2}{c}{\textbf{iter}} & \multicolumn{2}{c}{\textbf{sol. rank}} & \multicolumn{2}{c}{\textbf{runtime (s)}} \\
		& \ikpik & \lrpcg & \ikpik & \lrpcg & \ikpik & \lrpcg \\
		\midrule
		$\alpha=0.1$     &     2 &     1 &    40 &    31 &    218 &   197 \\
		$\alpha=0.2$	 &     2 &     1 &    40 &    31 &    235 &   203 \\
		$\alpha=0.25$	 &     3 &     2 &    60 &    40 &    343 &   537 \\
		\midrule
		$\ell=1$		 &     4 &     3 &    80 &    88 &    530 &   790 \\
		$\ell=5$		 &     3 &     2 &    60 &    50 &    347 &   387 \\
		$\ell=10$		 &     2 &     1 &    40 &    31 &    220 &   208 \\
		\midrule
		$\sig^2=10^{-2}$ &     2 &     1 &    40 &    31 &    212 &   219 \\
		$\sig^2=10^{-3}$ &     3 &     2 &    60 &    50 &    363 &   539 \\
		$\sig^2=10^{-4}$ &     4 &     3 &    80 &    79 &    532 &   771 \\
		\bottomrule
	\end{tabular}
\end{table}

\begin{table}[h]
	\centering
	\caption{Degree-Weighted Average}
	\label{tab:degree_weighted_average}
	\begin{tabular}{ c | r r  r r  r r}
		\toprule
		\textbf{degree-weighted} & \multicolumn{2}{c}{\textbf{iter}} & \multicolumn{2}{c}{\textbf{sol. rank}} & \multicolumn{2}{c}{\textbf{runtime (s)}} \\
		\textbf{average} & \ikpik & \lrpcg & \ikpik & \lrpcg & \ikpik & \lrpcg \\
		\midrule
		$\ell=1$		 &    3 &    2 &     60 &    30 &    112 &    213 \\
		$\ell=3$		 &    5 &    4 &    100 &    47 &    175 &    456 \\
		$\ell=5$		 &    5 &   24 &    100 &    55 &    144 &   3325 \\
		\midrule
		$\sig^2=10^{-2}$ &    3 &    2 &     60 &    43 &     97 &    242 \\
		$\sig^2=10^{-3}$ &    3 &    2 &     60 &    31 &     94 &    230 \\
		$\sig^2=10^{-4}$ &    2 &    1 &     40 &    22 &     76 &     93 \\
		\bottomrule
	\end{tabular}
\end{table}

The results show that the \lrpcg\ almost always requires fewer iteration steps than \ikpik. The rank of a solution obtained by \lrpcg\ is lower than that of a solution obtained by \ikpik. In return \ikpik\ has almost always a shorter runtime than \lrpcg. As we noticed, \lrpcg\ handles poorly conditioned systems better and converges more reliably. In the tested ranges of the hyperparameters $\alpha, \ell$ and $\sigma^2$, no significant differences in the results were observed.

For parameters leading to very badly conditioned systems, we observed instabilities in \ikpik, which seem to originate from the solvers for the small projected Sylvester equations. While we observe that both methods are feasible for models up to more than 28,000 nodes, the degree-weighted average can take over 3 million in most cases in 6 minutes on our machine.

\section{Conclusion}\label{sec:conclusion}
We demonstrated that, under the separability assumption of the covariance function of a MOGP and given low-rank data, the posterior mean can be efficiently computed by solving the Stein equation using the state-of-the-art \ikpik\ method. While \ikpik\ is an excellent solver, the proposed \lrpcg\ method--preconditioned with \ikpik--offers better control over the solution rank, albeit at the cost of slower performance in general. Additionally, under the strong but in some cases realistic assumption of stationary data, we showed that the degree-weighted average covariance leads to an efficient computation of the solution of the Stein equation.
The computation of the posterior variances and the hyperparameter optimization will be tasks of future work.

\section*{Acknowledgments}
S.E. thanks the ESF Plus of the European Union and the Free State of Saxony for funding and supporting this paper and the project.
S.E. used ChatGPT to polish the language of parts of this paper.

\newpage

\appendix

\section{Submatrices of filter matrices}\label{sec:submatrices}
We here want to consider how to construct the output covariance matrices in the general case with different training and target nodes $V\neq V^\ast$, where the matrix $\KKO$ is a submatrix of the filter matrix of the whole graph. This leads to a problem because we can not simply take a submatrix of, say $(\II+\alpha \LL)^2$, and solve by it since it is not the same as the submatrix of $(\II+\alpha \LL)^{-2}$. A possible approach is by \textit{selection matrices} which extract rows and columns of another matrix by multiplication. For extracting rows with indices $r_1, \dots, r_p$ we multiply $\RR$ from the left and for columns with indices $c_1, \dots, c_q$ we multiply $\CC$ from the right. Let $\ee_k$ be the $k$-th unit vector. We define 
\begin{align}
	\RR= \vek{\ee_{r_1}^\top \\ \vdots \\ \ee_{r_p}^\top}, \qquad
	\CC= \vek{\ee_{c_1} & \dots & \ee_{c_q}}.
\end{align}
We use these matrices to avoid the costly construction of the dense inverse matrix of the matrix filters.

\subsection*{Global filter}
While $\KKO$ is applied within the Sylvester equation we need to extract it as submatrix from
\begin{align}
	\BBT = (\II+\alpha\LL)^{-2}.
\end{align}
We achieve this by using the row and column selection with rows and columns from $\Vc$
\begin{align}
	\KKO = \RR (\II+\alpha\LL)^{-2}\CC.
\end{align}
In an initialization step we suppose to do a Cholesky decomposition of $\II+\alpha\LL$. For multiplication by a vector $\bb$ of fitting size we multiply first by $\CC$, solving by the Cholesky solver twice and multiply with $\RR$ at the end.

If we want to compute $\KKO^{-1}\bb$ we have to solve
\begin{align}\label{eqn:globalFilterInverse}
	\RR (\II+\alpha\LL)^{-2}\CC\xx = \bb.
\end{align}
We use that \Cref{eqn:globalFilterInverse} is the Schur complement of the following system
\begin{align}
	\begin{bmatrix}
		0 & -\RR \\
		\CC & \XX
	\end{bmatrix}\vek{\xx \\ \ast} = \vek{\bb \\ 0}.
\end{align}
We LU factorize the block matrix in the initialization step and extract $\xx$ from the solution vector in each solving step.

\subsection*{Local average filter}
For gaining the submatrix $\KKO$ from the total covariance matrix we use again selection matrices
\begin{align}
	\KKO = \RR(\II+\alpha\DD)^{-1}(\II+\alpha\WW)^2(\II+\alpha\DD)^{-1}\CC.
\end{align}
Computing $\FF=(\II+\alpha\WW)(\II+\alpha\DD)^{-1}$ maintains the sparsity and is done quickly because of the diagonal structure of $\II+\alpha\DD$.

For solving by $\KKO$ we use the Cholesky decomposition of $\RR\FF^\top\FF\CC$ and solve then by the Cholesky factor.

\end{document}